\documentclass[11pt]{article}
\usepackage{latexsym}
\usepackage{amsmath}

\renewcommand{\epsilon}{\varepsilon}

\usepackage{amsmath,bm}
\usepackage{amssymb}
\usepackage{amsthm}
\usepackage[ansinew]{inputenc}
\usepackage{graphicx}
\usepackage{epsfig}
\usepackage{dsfont}
\usepackage{bbm}
\usepackage{subfig}

\usepackage[authoryear]{natbib}

\newtheorem{satz}{Theorem}[section]

\newtheorem{lem}[satz]{Lemma}

\newtheorem{kor}[satz]{Corollary}
\newtheorem{rem}[satz]{Remark}
\newtheorem{assumption}[satz]{Assumption}

\def\3{\ss}

\def\ze{\mathbb{Z}}

\def \R{I \!\! R}   
\def \N{I \!\! N}   
\def \Z{Z \!\!\! Z} 

\newcommand{\E}{\mathbbm{E}}

\newcommand{\bea}{\begin{eqnarray*}}
\newcommand{\eea}{\end{eqnarray*}}
\newcommand{\be}{\begin{eqnarray}}
\newcommand{\ee}{\end{eqnarray}}

\newcommand{\var}{ \mbox{\sl Var} \ }

\newcommand{\ba}{\begin{array}}
\newcommand{\ea}{\end{array}}
\newcommand{\cum}{\text{\rm cum}}

\newcommand{\Cov}{\text{\rm Cov}}
\newcommand{\Var}{\text{\rm Var}}
\def\3{\ss}

\def\ze{\mathbb{Z}}

\allowdisplaybreaks

\begin{document}

\title{Measuring stationarity in long-memory processes}

\author{Kemal Sen, Philip Preu\ss, Holger Dette \\
Ruhr-Universit\"at Bochum \\
Fakult\"at f\"ur Mathematik \\
44780 Bochum \\
Germany \\
{\small email: kemal.sen@ruhr-uni-bochum.de}\\
{\small email: philip.preuss@ruhr-uni-bochum.de}\\
{\small email: holger.dette@ruhr-uni-bochum.de}
}

 \maketitle
\begin{abstract}
In this paper we consider the problem of measuring stationarity in locally stationary long-memory processes. 
We introduce an $L_2$-distance between the spectral density of the locally stationary process and its best approximation under the assumption of stationarity. The distance is estimated by a numerical approximation of the integrated spectral periodogram and asymptotic normality of the resulting estimate is established. 
The results can be used to construct a simple test for the hypothesis of stationarity in locally stationary long-range dependent processes. We also propose a bootstrap procedure to improve the approximation of the nominal level and prove its consistency. Throughout the paper, we will work with Riemann sums of a squared periodogram instead of integrals (as it is usually done in the literature) and as a by-product of independent interest it is demonstrated that the two approaches behave differently in the limit.
\end{abstract}



\section{Introduction}
\def\theequation{1.\arabic{equation}}
\setcounter{equation}{0}

\begin{figure}
\begin{center}
\includegraphics[width=0.3\textwidth]{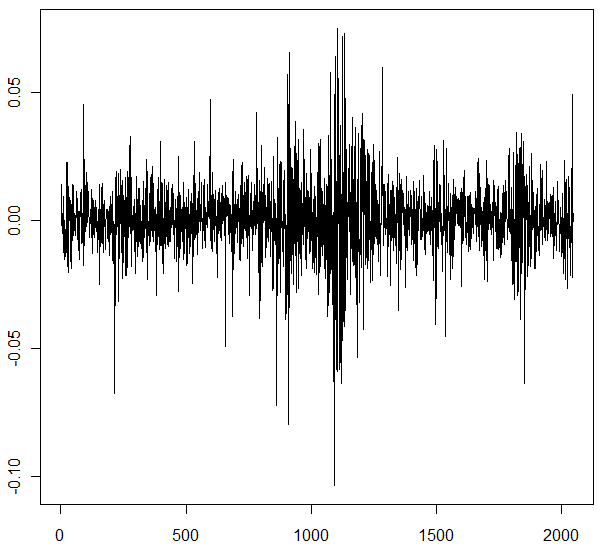}~~
\includegraphics[width=0.3\textwidth]{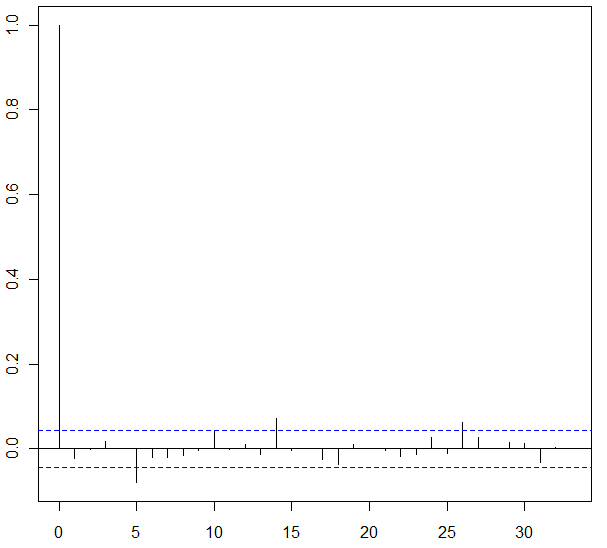}~~
\includegraphics[width=0.3\textwidth]{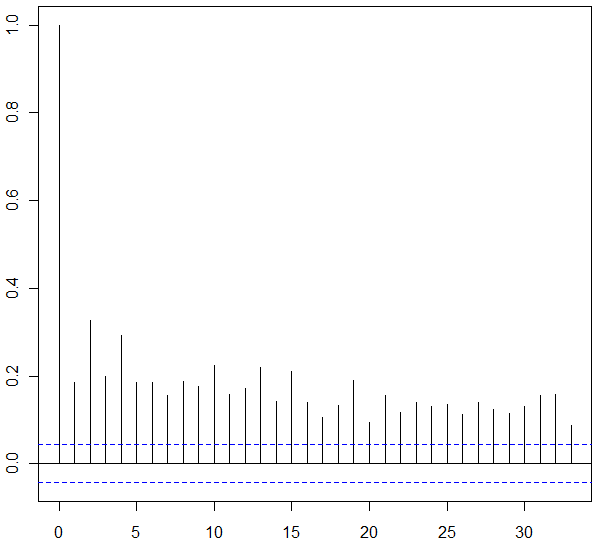}~~
  \caption{ \label{ibm} {\it Left panel: log-returns of the IBM stock between June 9th 2004 and  July 24th 2012, middle panel: ACF of the log-returns $X_t$, right panel: ACF of the squared returns $X_t^2$} }
\end{center}
\end{figure}

The assumption of (second-order) stationarity is quite common in the analysis of time series data like wind speeds, computer network traffic or stock returns. This condition allows for a well developed statistical analysis, and there exist numerous books and articles dealing with parameter estimation or forecasting techniques. However, under the assumption of stationarity many real world phenomena can only be described by complicated and less intuitive models. A typical example can be found in the left panel of Figure \ref{ibm} which shows $2048$ log-returns of the IBM stock between June 9th 2004 and  July 24th 2012. We observe that the autocovariance function (ACF)
 $\gamma(k)=\Cov(X_0,X_k)$ of the log-returns $X_t$ is converging to zero very ``fast'' as $k \rightarrow \infty$, while this is not the case for the ACF of the squared returns $X_t^2$  [see the middle and right panel in Figure \ref{ibm}]. The latter effect serves as the usual motivation to employ stationary long-memory models in the analysis of stock volatilities [see \cite{breicrali}]. This means that stationary processes satisfy
\be \label{asymptgamma}
\gamma(k) \sim C k^{2d-1}, \quad k \rightarrow \infty
\ee
for some $d \in (0,0.5)$,  which is called the long-memory parameter. Examples which fit into this framework  are the well-known FARIMA($p,d,q$)-models introduced by \cite{grangerjoyeux} and \cite{hosking}. However, these kinds of processes are not very intuitive and it was suggested by several authors that one should use simple but non-stationary ``short-memory'' models  instead [see for example \cite{mikosch2004}, \cite{starica2005}, \cite{fryzlewicz2006} or \cite{motivation2} among others]. Therefore an important  question of interest in this context is, if the data should be analyzed by a stationary long-range dependent model or by a non-stationary ``short-memory''  model.

In the present paper we propose a measure of stationarity in long-range dependent locally stationary processes, which is used for the construction of a consistent test for the hypothesis of stationarity.
Since the assumption of stationarity is crucial in the application of various statistical tools, there exist several procedures to validate this condition
in the context of short-memory processes. A first test for stationarity in locally stationary processes [as introduced by  \cite{dahlhaus1997}]
was proposed by \cite{sacneug2000} and  is  based on the estimation of wavelet coefficients by a localized version of the periodogram. 
 \cite{paparoditis2009,paparoditis2010} suggested an $L_2$-distance between the estimated spectral 
densities under the assumptions of stationarity and of local stationarity, and \cite{rao2010} developed a Portmanteau-type test statistic to detect deviations from stationarity. 
Besides the choice of a window 
width  for the localized periodogram, which is inherent in essentially any  statistical inference for locally stationary processes, all these methods require the choice of at least one additional smoothing parameter, such as the order of the wavelet expansion, a bandwidth for the estimation of the spectral density or the order in a Portmanteau-type test. 
\cite{detprevet2010} developed tests for stationarity in the framework of locally stationary processes which only require the choice of one regularization parameter, namely the   window length for the localized periodogram in the estimation procedure. These authors considered the $L_2$-distance
\be
\label{L2approach}
D^2:=\int_0^1\int_{-\pi}^\pi \left(f(u,\lambda)-\int_0^1 f(v,\lambda) dv \right)^2 d\lambda du
\ee

between the time varying spectral density $f(u,\lambda)$ and its approximation $\lambda \mapsto \int_0^1 f(v,\lambda) dv$ through a spectral density corresponding to a stationary process. 
It is easy to see that the process is stationary (i.e. the time varying spectral density does not depend on $u$) if and only if $D^2=0$, and $D^2$ can be considered as a measure of deviation from stationarity in the frequency domain.
 This quantity corresponds to the measure used in \cite{paparoditis2009}, but unlike to this author, \cite{detprevet2010} estimated $D^2$ directly via Riemann sums of the (squared) local periodogram instead of a smoothed local periodogram and thus avoided the choice of a smoothing parameter.  \cite{detprevet2011b} proposed an alternative measure for deviations from stationarity which is based on the Kolmogorov-Smirnov distance
\be
\label{KSapproach}
D_{KS}:=\sup_{(v,\omega ) \in [0,1]^2} \frac{1}{\pi} \Bigl | \int_0^v\int_{-\pi \omega}^{\pi \omega} f(u,\lambda)d \lambda d u-v\int_{-\pi \omega}^{\pi \omega} \int_0^1 f( u,\lambda) du d \lambda \Bigr |
\ee

[see also \cite{dahlhaus2009}]. Both approaches have their pros and cons. 
In particular tests based on the distance (\ref{KSapproach}) are $\sqrt{T}$-consistent (here $T$ denotes the sample size). On the other hand it is well known that 
- although such tests are consistent against alternatives converging to the null hypothesis at a parametric rate -
Kolmogorov-type and related tests greatly  weigh down contributions from high frequency components [see \cite{ghohuang1991}, \cite{euriccia1992} or \cite{fan1996}].
Moreover, the limiting distribution of Kolmogorov-Smirnov-type test statistics is usually not known. In principle this problem can be solved by bootstrap methods. However in many cases this yields to a loss of power. 
On the other hand, tests based on the $L_2$-approach can often use critical values from the standard normal distribution. 

As all procedures which have been suggested so far for discriminating between stationarity and non-stationarity, 
the tests proposed by \cite{detprevet2010} and \cite{detprevet2011b} are only applicable to short-memory processes, and the development of a corresponding methodology in the context of long-range dependence is missing. In fact, although stationary long-memory models are employed numerously in practice, there do not exist many tests for the hypothesis of stationarity which include these processes. \cite{horvarth2006}, \cite{sibbertsen2009} and \cite{aeneas} consider CUSUM and Wilcoxon type tests to discriminate between long-range dependence and one change in mean. However a change with respect to the mean is of course only the simplest possible deviation from stationarity and there is particular interest in measuring deviations in the dependency structure over time as well.

Recently, \cite{prevet2012} developed a test for stationarity which includes the long-range dependent case and is based on the distance \eqref{KSapproach}. As mentioned in the previous paragraph there exist several situations where this approach is not the best and for this reason we consider in this paper an alternative test which is based on the measure defined in \eqref{L2approach}. For this purpose, we estimate the integrals in the distance $D^2$ in (\ref{L2approach}) by Riemann sums where the unknown spectral densities are replaced by periodograms. For the resulting statistic we will show
that an appropriately standardized statistic converges to a standard normal distribution if the (possibly time varying) long-memory parameter $d(u)$ is smaller than $1/8$. These results are used to develop a bootstrap procedure for the approximation of the limit distribution and to prove its consistency in the general case.

Although that the proof of asymptotic normality seems to be more of theoretical nature, because
the bootstrap procedure derived in the second part of the paper can in principle also be applied in the case $d(u) < 1/8$, these results are of interest from several perspectives. Firstly, several arguments used in the proof of asymptotic normality are also required in the proof of consistency of the bootstrap procedure and easier illustrated in the unconditional case. Secondly, and most important, the estimate $\hat{D}_T^2$ of $D^2$ is based on estimates of the integrated and integrated squared spectral density $\int_{0}^{1} f(u,\lambda) \,du $ and $\int_{0}^{1}\int_{-\pi}^{\pi} f^2(u,\lambda) \,d\lambda \,du$, respectively. For this purpose we use Riemann sums of the squared periodogram instead of not computable integrals as it is usually done in the literature [see \cite{taniguchi1980}, \cite{foxtaqqu} and \cite{palole2010} among others]. Although one might expect that both estimators exhibit a similar behavior with respect to weak convergence, it is demonstrated in Section \ref{schwache Konvergenz} that this is not the case in the present context. A similar observation was also made by \cite{deochen2000} in the case of short-memory stationary processes. To the best of our knowledge, even in the (much simpler) stationary case, Riemann sums of a squared periodogram have not been considered in the literature for the long-range dependent case. 

The remaining part of this paper is organized as follows: In Section $2$, we introduce the necessary notation and define an empirical measure of stationarity. 
In Section \ref{schwache Konvergenz}, we prove that an appropriately standardized version of this measure converges weakly to a standard normal distribution if the time varying long-memory parameter is smaller than $1/8$. 
In Section $4$, we present a bootstrap procedure to approximate the distribution of $\hat{D}_T^2 $ and prove its consistency. The finite sample properties are investigated in Section $5$.
Finally, we defer all technical details to appendices in Section $6$ and $7$.

\section{Measuring stationarity in locally stationary long-memory processes}
\def\theequation{2.\arabic{equation}}
\setcounter{equation}{0}

In order to obtain a measure of stationarity including the long-range dependent case, we require a set-up which is flexible enough to cover stationary long-memory processes and a reasonable time-varying extension of it as well. For this reason, we  consider the following theoretical framework of a locally stationary long-memory process [similar approaches can be found in \cite{beran2009}, \cite{palole2010} and \cite{vonsachs2011}].
\\
\noindent
\begin{assumption}
\label{ass1}
Let $\{X_{t,T}\}_{t=1,...,T}$ denote a sequence of stochastic processes  which have a MA($\infty$) representation of the form
\begin{align}
\label{proc}
X_{t,T}=\sum_{l=0}^{\infty}\psi_{t,T,l}Z_{t-l}, \quad t=1,\ldots,T, 
\end{align}
such that
\be
\label{summesqendl}
\sup_{T \in \mathbb{N}}\sup_{t \in \{ 1, \ldots, T\}} \sum_{l=0}^\infty \psi_{t,T,l}^2<\infty,
\ee
where $\{Z_{t}\}_{t \in \mathbb{Z}} $ are independent and standard normal distributed random variables. We further assume the following conditions.
\begin{itemize}
\item[1)] 
There exist twice continuously differentiable functions $\psi_l: [0,1] \rightarrow \R$ ($l\in  \Z$)  such that
\be
\label{apprbed}
\sup_{t = 1, \ldots,T}\Bigl |\psi_{t,T,l}-\psi_l(t/T) \Bigr | \leq \frac{C}{T I(l)^{1-d_{\infty}}} \quad \forall l \in \N
\ee
and
\be
\psi_l(u)=\frac{a(u)}{I(l)^{1-d(u)}}+O\left(\frac{1}{I(l)^{2-d_{\infty}}}\right) \label{apprpsi}
\ee
holds uniformly in $u$ as $l \rightarrow \infty$, where $d,a :[0,1] \rightarrow \R_+$ are
twice differentiable functions, $C \in \R^+$ and $d_{\infty}=\sup_{u \in [0,1]}d(u)$ are constants and $I(x):=|x| \cdot 1_{\{x\not= 0\}} +1_{\{x=0\}}$.

\item[2)] The time varying spectral density 
\be
\label{tvspectraldensity}
f(u,\lambda):=\frac{1}{2\pi}\Bigl|\sum_{l=0}^\infty \psi_l(u)\exp(-i \lambda l)\Bigr|^2
\ee
is twice continuously differentiable on $(0,1)\times (0,\pi)$. Moreover, $f(u,\lambda)$ and all its partial derivatives up to order two are continuous on $[0,1]\times (0,\pi]$. 
\item[3)] There exists a constant $C \in \R^{+}$, which is independent of $u$ and $\lambda$, such that for $l \not= 0$
\be
\sup_{u \in (0,1)}|\psi_l'(u)| &\leq &  \frac{C\log |l|}{|l|^{1-d_{\infty}}},
\label{2.1b} \\
\sup_{u \in (0,1)}|\psi_l''(u)| &\leq &  \frac{C\log^2 |l|}{|l|^{1-d_{\infty}}}.
\label{2.1c}
 \ee
In addition, we assume
\be
\label{partialu}
\sup_{u \in (0,1)} \Big|\frac{\partial}{\partial u} f(u,\lambda)\Big| &\leq & \frac{C \log(\lambda)}{\lambda^{2d_{\infty}}}, \\
\label{partialu2}
\sup_{u \in (0,1)} \Big|\frac{\partial^2}{\partial u^2} f(u,\lambda)\Big| &\leq & \frac{C \log^2(\lambda)}{\lambda^{2d_{\infty}}}.
\ee
\end{itemize}
\end{assumption}
For the sake of a transparent notation, we will use $C \in \R^{+}$ as a universal constant throughout this paper. Note that the process is stationary if $\psi_{l,t,T}=\psi_l$ for all $l,t, T \in \N$. Condition \eqref{summesqendl} ensures that the infinite sum in \eqref{proc} exists in the $L_2$ sense, and \eqref{apprbed} means that the process $X_{t,T}$ can be approximated by a stationary model on a small time interval. It is also worthwhile to mention that the assumption of Gaussianity is only imposed here to simplify technical arguments [since they are quite involved even in this case]; see Remark \ref{Z_t non-gauss} for more details. 
Next, we consider the process
\be
\label{localappr}
X_t(u):=\sum_{l=0}^\infty \psi_l(u) Z_{t-l}
\ee
in order to visualize  some properties of a locally stationary long-memory process. Firstly, $X_t(u)$ is stationary for every fixed $u \in [0,1]$ and analogously to the stationary case, the condition (\ref{apprpsi}) implies the existence of bounded functions $y_i:[0,1] \rightarrow \mathbb{R}$ $(i=1,2)$ such that
\bea
|\Cov(X_t(u),X_{t+k}(u))| \sim \frac{y_1(u)}{k^{1-2d(u)}} \quad \text{as } k \rightarrow \infty 
\eea
and
\be
\label{apprpsifolgerung}
f(u,\lambda) \sim \frac{y_2(u)}{\lambda^{2d(u)}} \quad \text{ as } \lambda \rightarrow 0;
\ee
[see \cite{palole2010}].  Consequently, the autocovariance function $\gamma(u,k)=\Cov(X_0(u),X_k(u))$ is not absolutely summable and the time varying spectral density $f(u,\lambda)$ has a pole at $\lambda=0$ for any $u \in [0,1]$. 

As an example which fits in this theoretical set-up we consider the time varying FARIMA($p,d,q$) model defined by the equation
\be
\label{tvFARIMApq}
a(t/T,B)(1-B)^{d(t/T)} X_{t,T}=b(t/T,B)Z_t, \quad t=1,...,T,
\ee
where $B$ denotes the backshift operator,
\begin{eqnarray*}
a(u,z):=1-\sum_{j=1}^p a_j(u) z^j,\quad b(u,z):=1+\sum_{j=1}^q b_j(u) z^j
\end{eqnarray*}
 for given functions $a_j, b_j:[0,1] \rightarrow \mathbb{R}$, and $d:[0,1] \rightarrow (0,d_{\infty}]$ is twice continuously differentiable with $d_{\infty}<1/2$. It is shown in \cite{prevet2012} that under certain regularity conditions on the functions $a_j$, $b_j$, these equations have a locally stationary solution in the sense of Assumption \ref{ass1}. If the functions $a_j$, $b_j$ and $d$ do not depend on $u$, \eqref{tvFARIMApq} corresponds to the common FARIMA($p,d,q$) equation [see for example \cite{palmabuch} for conditions for the existence of a solution], which is included in our theoretical framework. \\
For the construction of an estimate of the quantity \eqref{L2approach} we note that
\begin{eqnarray}
D^2=2\pi F_1-4\pi F_2,
\label{D2alt}
\end{eqnarray}
where
\be \label{2.4a}
F_1 &:=&
 {\frac{1}{2 \pi}}
\int_0^1 \int_{-\pi}^{\pi}  f^2(u,\lambda)  d \lambda du, \\
 \label{2.4b}
 F_2 &:=& \frac{1}{4\pi}\int_{-\pi}^{\pi} \Bigl( \int_0^1 f(u,\lambda) du \Bigr)^2 d \lambda.
\ee
Consequently, it follows from $\eqref{apprpsifolgerung}$ that the distance $D^2$ is only well defined if $d_{\infty} < \frac{1}{4}$. We assume without loss of generality that the sample size $T$ can be decomposed into $M$ blocks with length $N$ (i.e. $T=NM$) where $N$ and $M$ are positive integers and $N$ is even. A rough estimator for the time varying spectral density $f(u,\lambda)$ is then given by the local periodogram at the rescaled time point $u \in [0,1]$ which is defined by
\begin{align*}
 I_N(u,\lambda):=\frac{1}{2 \pi N} \Bigl|\sum_{s=0}^{N-1} X_{\lfloor uT \rfloor-N/2+1+s,T} \exp(-i \lambda s) \Bigr|^2,
\end{align*}
where $X_{j,T} =0$ for $j \not \in \{1,\ldots , T\}$ [see \cite{dahlhaus1997}]. This is the usual periodogram computed from the observations  $X_{\lfloor uT \rfloor-N/2+1,T},\ldots, X_{\lfloor uT \rfloor+N/2,T}$, and it can be shown that it is asymptotically unbiased for the time-varying spectral density $f(u,\lambda)$ if $N \rightarrow \infty$ and $N=o(T)$. However, $ I_N(u,\lambda)$ is not consistent just as the usual periodogram. In addition, $I_N(u,\lambda)^2$ is an unbiased (but not consistent) estimate of the quantity $2f^2(u,\lambda)$ instead of $f^2(u,\lambda)$.

We now construct empirical versions of \eqref{2.4a} and \eqref{2.4b} by replacing the integrals through appropriate Riemann-sums and substitute $f(u,\lambda)$ and $f^2(u,\lambda)$ by $ I_N(u,\lambda)$ and $ I_N(u,\lambda)^2/2$, respectively. For this purpose, we define the rescaled mid-points of the $M$ blocks 
\begin{eqnarray*}
u_j:=\frac{t_j}{T}  := \frac{N(j-1)+N/2}{T} \quad (j=1,\ldots,M)
\end{eqnarray*}
 and consider the statistics
\be
\label{2.6}
\hat F_{1,T} &:=& \frac{1}{T} \sum_{j=1}^M \sum_{k=1}^{\lfloor \frac{N}{2}\rfloor}   I_N(u_j,\lambda_{k,N})^2, \\
\hat F_{2,T} &:=& \frac{1}{N}\sum_{k=1}^{\lfloor \frac{N}{2} \rfloor} \Bigl(\frac{1}{M}\sum_{j=1}^M I_N(u_j,\lambda_{k,N})\Bigr)^2,
\label{2.7}
\ee
where $\lambda_{k,N} = 2 \pi k/N$ denote the usual Fourier frequencies. The empirical measure of stationarity \eqref{L2approach} is finally given by
\be
\label{2.8}
\hat D_T^2 := 2\pi \hat F_{1,T} - 4\pi \hat F_{2,T}.
\ee

We would like to point out here that it is far from obvious that $\hat D_T^2$ is a consistent estimator of $D^2$. In general it is  not true that an integrated function of the periodogram converges to the corresponding integrated function of the spectral density. This - at a first glance - is a counterintuitive property of the integrated periodogram and was already observed by 
 \cite{taniguchi1980} in the context of stationary short-memory processes. These problems are also visible here as well as we require a multiple of $I_N(u,\lambda)^2$  to obtain an asymptotically unbiased estimator for $f^2(u,\lambda)$. In the following section we will prove consistency of $\hat D_T^2 $ and study its weak convergence.

\section{Consistency and weak convergence}
\label{schwache Konvergenz}
\def\theequation{3.\arabic{equation}}
\setcounter{equation}{0}
Throughout this paper, the symbols $\xrightarrow{\text{ } \mathcal{P}\text{ }}$ and $  \xrightarrow{\text{ } \mathcal{D}\text{ }}$ denote convergence in probability and weak convergence, respectively. In order to specify the bias of $\hat{F}_{1,T}$ and $\hat{F}_{2,T}$ we define
\bea
F_{1,T} := \frac{1}{2\pi M}  \sum_{j=1}^M \int_{-\pi}^{\pi}   f^2(u_j,\lambda) \,d\lambda, \quad F_{2,T} := \frac{1}{4\pi}\int_{-\pi}^{\pi}  \Bigl(\frac{1}{M}\sum_{j=1}^M f(u_j,\lambda)\Bigr)^2  \,d\lambda,
\eea

and obtain the following results.
\begin{satz}\label{thm0}
Suppose Assumption \ref{ass1} holds with $\sup_{u \in [0,1]} d(u)<1/4$ and that the conditions
\be
N \rightarrow \infty,\quad \frac{N}{T} \rightarrow 0
\ee
are satisfied. Then 
$\hat{F}_{1,T}  \xrightarrow{\text{ }   \mathcal{P}\text{ }} F_{1},\hspace{.1cm} \hat{F}_{2,T}  \xrightarrow{\text{ }   \mathcal{P}\text{ }} F_{2} 
$
and in particular
\begin{eqnarray*}
\hat{D}^2_T \xrightarrow{\text{ } \mathcal{P}\text{ }} D^2
\end{eqnarray*}
as $T \rightarrow \infty $.
\end{satz} 
\begin{satz} \label{thm1} Suppose Assumption \ref{ass1} holds with $d_{\infty}=\sup_{u \in [0,1]} d(u)<1/8$ and that the conditions
\be
\label{assNM} 
N \rightarrow \infty,\quad \frac{N}{T} \rightarrow 0, \quad \frac{\sqrt{T}}{N^{1-4d_{\infty}}} \rightarrow 0
\ee
are satisfied. Then as $T \rightarrow \infty$ we have
$$
\sqrt{T} \bigl\{ ( \hat F_{1,T} , \hat F_{2,T}) ^T   - (F_{1,T}, F_{2,T} + d_{N,T})^T -\boldsymbol{C_T} \bigr \}
 \xrightarrow{\text{ }   \mathcal{D}\text{ }} \mathcal{N}(0,\boldsymbol{\Sigma}),
$$
 where the covariance matrix $\boldsymbol{\Sigma}$ and the constant $ d_{N,T}$ are given by
\be\label{3.3}  \boldsymbol{\Sigma}&=&\hspace{-.2cm} \begin{pmatrix} \frac{5}{ \pi} \int_{-\pi}^{\pi} \int_{0}^{1}f^4(u,\lambda) du d \lambda  & \frac{2}{\pi} \int_{-\pi}^{\pi} \Bigl( \int_0^1 f(u,\lambda)du \int_0^1 f^3(u,\lambda)du \Bigr) d \lambda \\ \frac{2}{\pi} \int_{-\pi}^{\pi} \Bigl( \int_0^1 f(u,\lambda)du \int_0^1 f^3(u,\lambda)du \Bigr) d \lambda & \frac{1}{\pi}\int_{-\pi}^{\pi} \Bigl( \Bigl( \int_0^1 f(u,\lambda)du \Bigr)^2 \int_0^1 f^2(u,\lambda)du \Bigr) d \lambda \end{pmatrix} \\
\label{3.3a}
 d_{N,T}&=&\frac{ 1}{  4\pi  M^2} \sum_{j=1}^M \int_{-\pi}^{\pi} f^2(u_j,\lambda) \,d\lambda,
 \ee
respectively, and the vector $\boldsymbol{C_T} \in \R^2$  is of order $O\left(N^2/T^2+\log(N)/(MN^{1-4d_{\infty}})\right)$. In particular, this term vanishes if the functions $\psi_l(u)$ are independent of $u$ for all $l \in \Z$ [i.e. the spectral density $f(u,\lambda)$ of the underlying process $X_{t,T}$ is independent of $u$].
 \end{satz}
A similar result for the short-memory situation has been derived by \cite{detprevet2010}. In contrast to their result, there appears an additional bias term $\boldsymbol{C_T}$ in Theorem \ref{thm1}. This term is negligible under the additional condition $N^2/T^{3/2}\rightarrow 0$  which holds under the stronger restriction $d_{\infty}<1/12$ due to \eqref{assNM}. 
On the other hand, under the null hypothesis of a time independent spectral density 
\begin{align} \label{H0}
H_0: f(u,\lambda) \text{ is independent of } u,
\end{align}


we have that $\boldsymbol{C_T}=0$ (this follows from the proof of Theorem \ref{thm1} in the Appendix). Since the covariance matrix \eqref{3.3} contains the integrated fourth power of the spectral density, we obtain from \eqref{apprpsifolgerung} that  Theorem \ref{thm1} is not valid whenever $d_{\infty} \geq 1/8$.
Writing $(C_{1,T},C_{2,T})^T:=\boldsymbol{C_T}$, a straightforward application of the Delta-method yields the following result.
\begin{kor}
\label{asymp normalvert}
Under the assumptions of Theorem \ref{thm1}, it holds
\be
\label{zentrb2}
\sqrt{T} \Bigl( \hat D_T^2  -   D_T^2 +  4\pi d_{N,T}+4 \pi C_{2,T}-2\pi C_{1,T} \Bigr) \xrightarrow{\text{ } \mathcal{D}\text{ }} \mathcal{N}(0,  \tau^2),
\ee
where $D_T^2 := 2\pi F_{1,T}-4\pi F_{2,T}$ and the asymptotic variance is given by
\be 
 \label{3.17}    \tau^2 &:=& 20\pi \int_{-\pi}^{\pi} \int_{0}^{1} f^4(u,\lambda) du d \lambda -32\pi \int_{-\pi}^{\pi} \Bigl( \int_0^1 f(u,\lambda)du \int_0^1 f^3(u,\lambda)du \Bigr) d \lambda \\
\nonumber && +16\pi\int_{-\pi}^{\pi} \Bigl( \Bigl( \int_0^1 f(u,\lambda)du \Bigr)^2 \int_0^1 f^2(u,\lambda)du \Bigr) d \lambda.
\ee
\end{kor}
Under the null hypothesis (\ref{H0}) we have $D_T^2=C_{1,T}=C_{2,T}=0$ and the asymptotic variance in \eqref{3.17} reduces to $\tau^2_{H_0} := 4 \pi \int_{-\pi}^{\pi} f^4 (\lambda)d \lambda$. The asymptotic bias $4 \pi d_{N,T}=\frac{ 2 \pi N}{T} F_{1,T}$ can easily be estimated by the statistic $\hat{B}_T:=\frac{2\pi N}{T} \hat F_{1,T}$  and we infer from Theorem \ref{thm1}
 $$
 \sqrt{T} \Bigl( \hat{B}_T - 4 \pi d_{N,T} \Bigr) = \frac{ 2\pi N}{T} \sqrt{T} \Bigl(\hat F_{1,T}- F_{1,T} \Bigr) \xrightarrow{\text{ } \mathcal{P}\text{ }} 0 .
 $$
Thus Slutzky's Lemma together with \eqref{zentrb2} yields
\be
\label{zentrnullhilf2}
 \sqrt{T} \Bigl( \hat D_T^2  + \hat{B}_T \Bigr) \xrightarrow{\text{ } \mathcal{D}\text{ }} \mathcal{N}(0,  \tau^2_{H_0})
\ee

under the null hypothesis. To construct an asymptotic level $\alpha$-test for stationarity, it therefore remains to estimate the variance $\tau^2_{H_0}$ in (\ref{zentrnullhilf2}), and an estimator for this quantity is given by $\hat \tau^2_{H_0} := 4  \pi^2 \hat \tau_1^2$ with 
\bea
 \hat \tau_1^2 &:=& \frac{1}{6T} \sum_{k=1}^{\lfloor \frac{N}{2} \rfloor} \sum_{j=1}^M I_N(u_j,\lambda_{k,N})^4.
\eea
The consistency of this estimator follows from the next theorem.
\begin{satz} \label{thm3} If the assumptions of Theorem \ref{thm1} are satisfied, we have
 \bea
 \hat \tau_1^2 & \xrightarrow{\text{ } \mathcal{P}\text{ }} &\frac{1}{\pi}  \int_{-\pi}^{\pi} \int_0^1 f^4(u,\lambda) du d \lambda.
\eea
\end{satz}
Combining \eqref{zentrnullhilf2} with Theorem \ref{thm3} yields
\be
\label{zentrnull}
 \sqrt{T} \Bigl( \hat D_T^2  + \hat{B}_T \Bigr)\Big/\sqrt{\hat \tau^2_{H_0}} \xrightarrow{\text{ } \mathcal{D}\text{ }} \mathcal{N}(0,1)
\ee

and therefore an asymptotic level $\alpha$-test is obtained by rejecting the null hypothesis  \eqref{H0} whenever
\begin{eqnarray}
 \sqrt{T} \Bigl( \hat D_T^2  + \hat{B}_T \Bigr) \Big/\sqrt{\hat \tau^2_{H_0}} \geq u_{1-\alpha} 
 \label{rule},
\end{eqnarray}
where $u_{1-\alpha}$ denotes the $(1-\alpha)$-quantile of the standard normal distribution. It follows from Theorem \ref{thm1} that this test is consistent, because the left hand side of (\ref{rule}) converges to infinity, whenever there exists a $\tilde \lambda \in [-\pi,\pi]$ such that the function $u \mapsto f(u,\tilde \lambda)$ is not constant.

\begin{rem}
\label{Z_t non-gauss}
{ \rm
If the innovation process $(Z_t)_{t \in \mathbb{Z}}$ in \eqref{proc} is not Gaussian, it can be shown that Corollary \ref{asymp normalvert} is still valid where the asymptotic variance $\tau^2$ in (\ref{3.17}) has to be replaced by
\begin{eqnarray*}
\tau_g^2&=& \tau^2+\frac{\kappa_4}{\kappa_2^2} \bigg \{  4 \int_{0}^{1} \bigg(\int_{-\pi}^{\pi} f^2(u,\lambda) \,d\lambda\bigg)^2 \,du+ 4 \int_{0}^{1} \bigg( \int_{-\pi}^{\pi}  f(u,\lambda) \bigg(\int_{0}^{1} f(\nu,\lambda) \,d\nu \ \bigg) \,d\lambda\bigg)^2 \,du \\
&& \hspace{1.5cm}-8\int_{0}^{1} \bigg( \int_{-\pi}^{\pi} f^2(u,\lambda) \,d\lambda \int_{-\pi}^{\pi} f(u,\lambda) \bigg( \int_{0}^{1} f(\nu,\lambda) \,d\nu \bigg)  \,d\lambda\bigg) \,du \bigg\}
\end{eqnarray*}
and $\kappa_2$ and $\kappa_4$ denote the second and fourth cumulants of the innovation process, respectively. In particular, under the null hypothesis of stationarity, it follows that $\tau_g^2=\tau^2= \tau^2_{H_0}$ and hence no adjustments in the asymptotic level $\alpha$-test in \eqref{rule} are necessary to address non normal distributed innovations.
}
\end{rem}

\begin{rem} \label{riemannstattintegral}  
{\rm
We note that for  locally stationary long-range dependent models the asymptotic variances of the statistics 
\bea
\tilde{F}_{1,T}=\frac{1}{4\pi M }\sum_{j=1}^{M}\int_{-\pi}^{\pi} I_N(u_j,\lambda)^2 \,d\lambda 
\eea
 and of $\hat F_{1,T}$, defined in (\ref{2.6}), are different. In fact it follows by similar arguments as given in the appendix that 
 \begin{eqnarray*}
\lim_{T \rightarrow \infty} T\hspace{0.05cm}\var\hspace{-0.1cm} (\tilde{F}_{1,T})=\frac{14}{3\pi}\int_{-\pi}^\pi\int_0^1 f^4(u,\lambda) du d\lambda,
 \end{eqnarray*}
  while  
 \begin{eqnarray*}
\lim_{T \rightarrow \infty} T\hspace{0.05cm}\var \hspace{-0.1cm}(\hat{F}_{1,T})=  \frac{5}{ \pi} \int_{-\pi}^{\pi} \int_{0}^{1}f^4(u,\lambda) du d \lambda
 \end{eqnarray*}
  by Theorem \ref{thm1}.
Moreover, similar arguments as given in the proof of this statement show that even in the stationary case 
 the asymptotic variance of the statistic $\int_{-\pi}^{\pi}   I_T(\lambda)^2 \,d\lambda$ and its discretized version $(2\pi/T)\sum_{k=1}^{T}  I_T(\lambda_{k,T})^2$ are not the same (here $ I_T(\lambda)$ denotes the usual periodogram and $\lambda_{k,T}=2\pi k/T$ are the Fourier frequencies). \cite{deochen2000} observed the same effect in the context of  stationary short-memory processes.


}
\end{rem}

\section{Critical values by resampling}
\def\theequation{4.\arabic{equation}}
\setcounter{equation}{0}

We now consider the more general set-up with $d_{\infty} < \frac{1}{4}$ as specified in Assumption \ref{ass1}. We will show that in this case a bootstrap procedure can be used to approximate the distribution of $\hat D_T^2$ under the null hypothesis \eqref{H0}. We employ the FARI($\infty$) bootstrap which was recently introduced  by \cite{prevet2012} and fits an FARIMA($p,\underline d,0$)-model to the data, where $p=p(T)$ converges to infinity with increasing sample size $T$. To prove consistency of this procedure, we require the following technical assumptions.

\begin{assumption}
\label{ass2}
For the stationary process $\{ X_t \}_{t \in \mathbb{Z}}$ with strictly positive spectral density $\lambda \mapsto \int_0^1 f(u,\lambda) du$, there exists a constant $ \underline d \in (0,1/4)$  such that the process
\be
\label{gl4.1}
Y_t=(1-B)^{\underline d} X_t
\ee
has an AR($\infty$)-representation of the form
\be
\label{statproc}
Y_t=\sum_{j=1}^\infty a_j Y_{t-j}+Z_t^{AR},
\ee
where  $\{Z_j^{AR}\}_{j\in \ze}$ denotes a Gaussian White Noise process with variance $\sigma^2>0$ and the coefficients in the representation $\eqref{statproc}$ satisfy
\be
\label{gl4.2}
\sum_{j=1}^\infty |a_j||j|^7<\infty, \\
1-\sum_{j=1}^\infty a_jz^j \not= 0 &\text{ for }& |z|\leq 1.
\ee
\end{assumption}
Note that under the null hypothesis of a time independent spectral density, it follows that $\underline d=d_{\infty}=d(u)$ for all $u \in [0,1]$, but under the alternative we usually have $\underline d \not= d_{\infty}$. The FARI($\infty$) bootstrap incorporates the following steps: 
First we choose a $p=p(T) \in \N$ to construct an estimator, say $\underline{\hat d}$, of the long-range dependence parameter $\underline d$ in model (\ref{gl4.1}). Secondly we calculate an estimator of
\be
\label{ar}
(a_{1,p},...,a_{p,p})= \underset{b_{1,p},...,b_{p,p}}{\operatorname{argmin}} \E \Bigl(Y_t-\sum_{j=1}^p b_{j,p}Y_{t-j} \Bigr)^2 ,
\ee
by fitting an AR($p$)-model to the data. In order to describe the main idea of our procedure in more detail, we introduce the ``true''  approximating process $Y_{t}^{AR}(p)$ by
\begin{align}
\label{apprprocess}
Y_{t}^{AR}(p)=\sum_{j=1}^p a_{j,p}Y_{t-j}^{AR}(p)+Z_{t}^{AR},
\end{align}
where the parameters $a_{j,p}$ are defined in (\ref{ar}) and $\{Z_{t}^{AR}\}_{t \in \mathbb{Z}}$ is a Gaussian White Noise process with mean zero and variance $\sigma_p^2=\E(Y_t-\sum_{j=1}^p a_{j,p}Y_{t-j})^2$. 
If $p=p(T) \rightarrow \infty$ the process $Y_{t}^{AR}(p)$ approximates $Y_t$ and therefore $(1-B)^{- {\underline d}} Y_{t}^{AR}(p)$ is ``close'' to the stationary process $X_t$ whose spectral density is given by $\lambda \mapsto \int_0^1 f(u,\lambda) du$. Under the null hypothesis of stationarity, this function coincides with the spectral density of $\{X_{t,T} \}_{t=1,\ldots, T}$. Hence, observing the data $X_{1,T},...,X_{T,T}$, the FARI($\infty$) bootstrap precisely works as follows: 
\begin{itemize}
\item[1)] Choose $p=p(T) \in \N$ and calculate $\hat \theta_{T,p}=(\underline{\hat d}, \hat \sigma_p^2, \hat a_{1,p},...,\hat a_{p,p})$ as the minimizer of 
\bea
 \frac{1}{T} \sum_{k=1}^{T/2} \left( \log f_{\theta_p}(\lambda_{k,T})+\frac{I_T(\lambda_{k,T})}{f_{\theta_p}(\lambda_{k,T})} \right),
\eea

where $\theta_p=(\underline{ d},  \sigma_p^2, a_{1,p},...,a_{p,p}),
$\bea
I_T(\lambda)=\frac{1}{2\pi T} \Bigl | \sum_{t=1}^T X_{t,T} \exp(-i \lambda t) \Bigr |^2
\eea
is the usual periodogram, and 
\bea
f_{\theta_p}(\lambda)= \frac{|1-\exp(-i \lambda)|^{-2 \underline d}}{2\pi} \times \frac{\sigma_p^2}{|1-\sum_{j=1}^p a_{j,p}\exp(-i \lambda j) |^2}
\eea

is the spectral density of a stationary FARIMA($p,\underline d,0$)-model. Note that the estimator $\hat \theta_{T,p}$ is the classical Whittle estimator of a stationary process [see \cite{whittle1}]. 
\item[2)] Calculate $Y_{t,T}=(1-B)^{\underline{\hat d}} X_{t,T}$ for $t=1,...,T$ and simulate a pseudo-series $Y_{1,T}^*,...,Y_{T,T}^*$ according to the model
\bea
Y_{t,T}^*&=&Y_{t,T}; \quad t=1,...,p ,\\
Y_{t,T}^*&=&\sum_{j=1}^p \hat a_{j,p}Y_{t-j,T}^*+\hat \sigma_p Z_j^*, \quad   p<t \leq T,
\eea
where $Z_j^*$ denotes an independent sequence of standard normal distributed random variables. 

\item[3)] Create the pseudo-series $X_{1,T}^*,...,X_{T,T}^*$ from the equation
\begin{eqnarray}
X_{i,T}^*=(1-B)^{-\underline{\hat d}} Y_{i,T}^*
\label{bootproc}
\end{eqnarray}
 and compute $\hat D_T^{2,*}$ in the same way as $\hat D_T^2$ where the original observations $X_{1,T},...,X_{T,T}$ are replaced by the bootstrap replicates $X_{1,T}^*,...,X_{T,T}^*$.
\end{itemize}

Our main theorem in this section describes the theoretical properties of this procedure.
\begin{satz}
\label{thmstat2}
Assume that the null hypothesis \eqref{H0} holds and let Assumption \ref{ass1} and \ref{ass2} be fulfilled. Furthermore, suppose that the conditions
\be
\label{assNM2} 
N \rightarrow \infty,\quad \frac{N}{T} \rightarrow 0, \quad \frac{T}{N^{1+\delta}} \rightarrow 0
\ee
are satisfied for some $0<\delta<1/2$, and assume for the growth rate (rate of convergence) of $p=p(T)$ the following:
\begin{itemize}
\item[i)]  
There exist sequences $p_{max}(T) \geq p_{min}(T) \xrightarrow{\text{ }T \rightarrow \infty \text{ }} \infty$ such that $p(T)\in [p_{min}(T),p_{max}(T)],$
\begin{eqnarray}
p_{max}^9(T) \log(T)^3N^\delta T^{-1}  =O(1), \label{asspNM} \\
\sqrt{T}p_{min}^{-9}(T)/\sqrt{\log(T)}=o(1).\label{asspNM2} 
\end{eqnarray}

\item[ii)] The condition
\be
\label{maxthetarate}
||\hat \theta_{T,p}- \theta_{p}||_\infty =O_P\left(\sqrt{\frac{\log(T)}{T}}\right)
\ee
is fulfilled uniformly with respect to $p$, where $\hat \theta_{T,p}$ denotes the estimator used  in step 1) of the bootstrap procedure and $\theta_{p}=(\underline{d}_p,\sigma_p^2,a_{1,p},...,a_{p,p})$ are the corresponding ``true'' parameters.
\end{itemize}
Then there exist  random variables $\hat D_{T,a}^{2}$ and $\hat D_{T,a}^{2,*}$ such that
\bea
&(a)\quad & \hat D_{T,a}^{2} \stackrel{\mathcal{D}}= \hat D_{T,a}^{2,*}, \\
&(b) \quad &  \Var(\hat D_{T}^{2})^{-1/2}\Big( \hat D_{T}^{2}-\hat D_{T,a}^{2} \Big)=o_P(1), \\
&(c) \quad &  \Var(\hat D_{T}^{2,*})^{-1/2}\Big( \hat D_{T}^{2,*}-\hat D_{T,a}^{2,*} \Big)=o_P(1), \\
&(d)\quad & \E \big|\hat D_{T,a}^{2,*}\big| =O\left(\frac{N^{\max(4\underline d-1/2,0)}}{\sqrt{T}}+\frac{\sqrt{\log(N)}1_{\{\underline d = 1/8\}}}{\sqrt{T}}+\frac{1}{N^{1-4 \underline d}}\right).
\eea
The estimate in (d) also holds if the null hypothesis \eqref{H0} is not satisfied.
\end{satz}

Note that conditions like \eqref{asspNM}-\eqref{maxthetarate} are standard in the context of parametric bootstraps [see for example \cite{bergpappolitis2010} or \cite{kreisspappol2011}] and a detailed discussion of them is given in \cite{prevet2012}. We now obtain an asymptotic level $\alpha$-test based on $\hat D_T^2$ as follows: Calculate B bootstrap replicates $\hat D_{T}^{2,*}$, denote by $(\hat D_T^{2,*})_{T,1},...,(\hat D_T^{2,*})_{T,B}$ the resulting order statistic and reject the null hypothesis whenever
\be
\label{test}
\hat D_T^2 > (\hat D_T^{2,*})_{T,\lfloor (1-\alpha)B \rfloor} .
\ee

Theorem \ref{thmstat2} and the argumentation in \cite{paparoditis2010} indicate that this procedure is valid for obtaining an asymptotic level $\alpha$-test. In order to prove this more formally, we follow \cite{bickel} by considering the Mallow metric $d_2(F,G)=\inf \sqrt{\E(X-Y)^2}$ between two distributions $F$ and $G$, where the infimum is taken over all pairs $(X,Y)$ of random variables with marginal distributions $F$ and $G$. Theorem \ref{thmstat2} then yields the following result which states that the test \eqref{test} has, in fact, asymptotic level $\alpha$.

\begin{satz} \label{mallowmetric}
Suppose the null hypothesis \eqref{H0} and the assumptions of  Theorem \ref{thmstat2} are satisfied. Then, as $T\rightarrow \infty$, the Mallow distance $d_2$ between the distributions of the random variables  
\begin{eqnarray*}
\hat D_T^2/\sqrt{\Var(\hat D_{T}^{2})} \mbox{ and }\hat D_T^{2,*}/\sqrt{\Var(\hat D_{T}^{2,*})}
\end{eqnarray*}
 converges to zero in probability.
\end{satz}

Consistency under the alternative follows since  Theorem \ref{thmstat2} d) yields that each bootstrap statistic $\hat D_T^{2,*}$ converges to zero while $\hat D_T^2$ exceeds some positive constant (for $T$ sufficiently large) due to Theorem \ref{thm0}.

\section{Finite sample properties}
\label{Simulationen}
\def\theequation{5.\arabic{equation}}
\setcounter{equation}{0}

In this section we examine the finite sample properties of the proposed decision rule \eqref{test}. An important problem is the choice of the window length $N$ 
for the calculation of the local periodogram and the  choice of the AR parameter $p$ in the bootstrap procedure. Throughout this section we choose $p$ as the minimizer of  the AIC criterion [see \cite{akaike1973}], which is defined by
\begin{align*}
\hat p = \text{argmin}_p \frac{2\pi}{T} \sum_{k=1}^{{T}/{2}} \Big( \log f_{\hat \theta(p)}(\lambda_{k,T})+\frac{I_T(\lambda_{k,T})}{f_{\hat \theta(p)}(\lambda_{k,T})} \Big) + \frac{p}{T}
\end{align*}
in the context of stationary processes due to \cite{whittle1} [here $f_{\hat \theta(p)}$ is the spectral density of the fitted stationary FARIMA($p,d,0$) process and $I_T$ is the usual stationary periodogram]. We therefore restrict ourselves to an analysis of the sensitivity  with respect to $N$ in the following, and it will turn out that the test (\ref{test}) using the FARI($\infty$) bootstrap exhibits a remarkable robustness with respect to the choice of $N$. All reported results of this section are based on $1000$ simulation runs and $200$ bootstrap replications.

\subsection{Size of the test}

\begin{figure}
\begin{center}
\includegraphics[width=1\textwidth,height=0.5\textwidth]{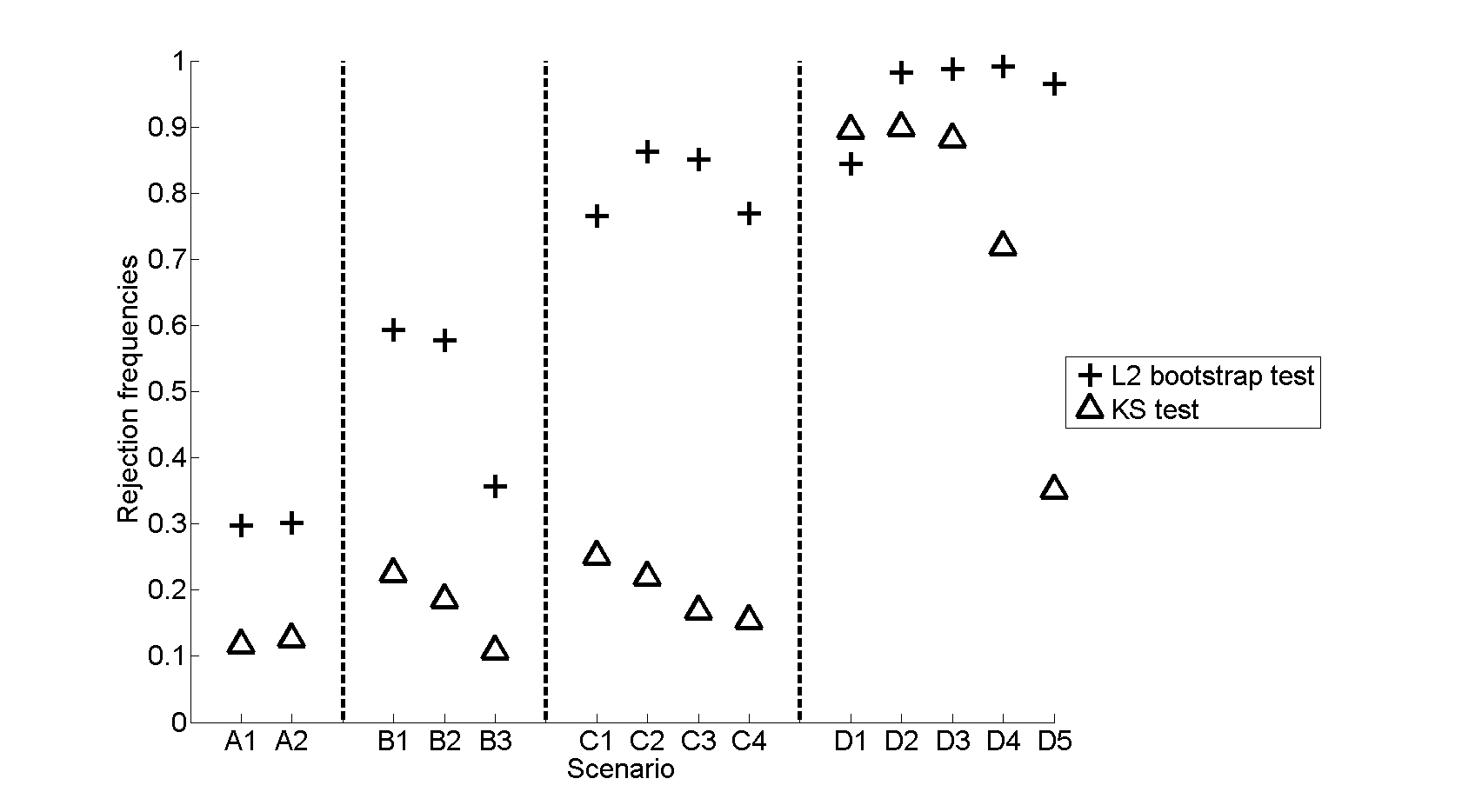}~~ \\
  \caption{ \label{alternative1} {\it Power of the test (\ref{test}) and the Kolmogorov-Smirnov test of \cite{prevet2012} at  5\% level for the model \eqref{alt1} under the scenarios A1-D5 from Table \ref{TableARd=0.1 mitboot}.} }
\end{center}
\begin{center}
\includegraphics[width=1\textwidth,height=0.5\textwidth]{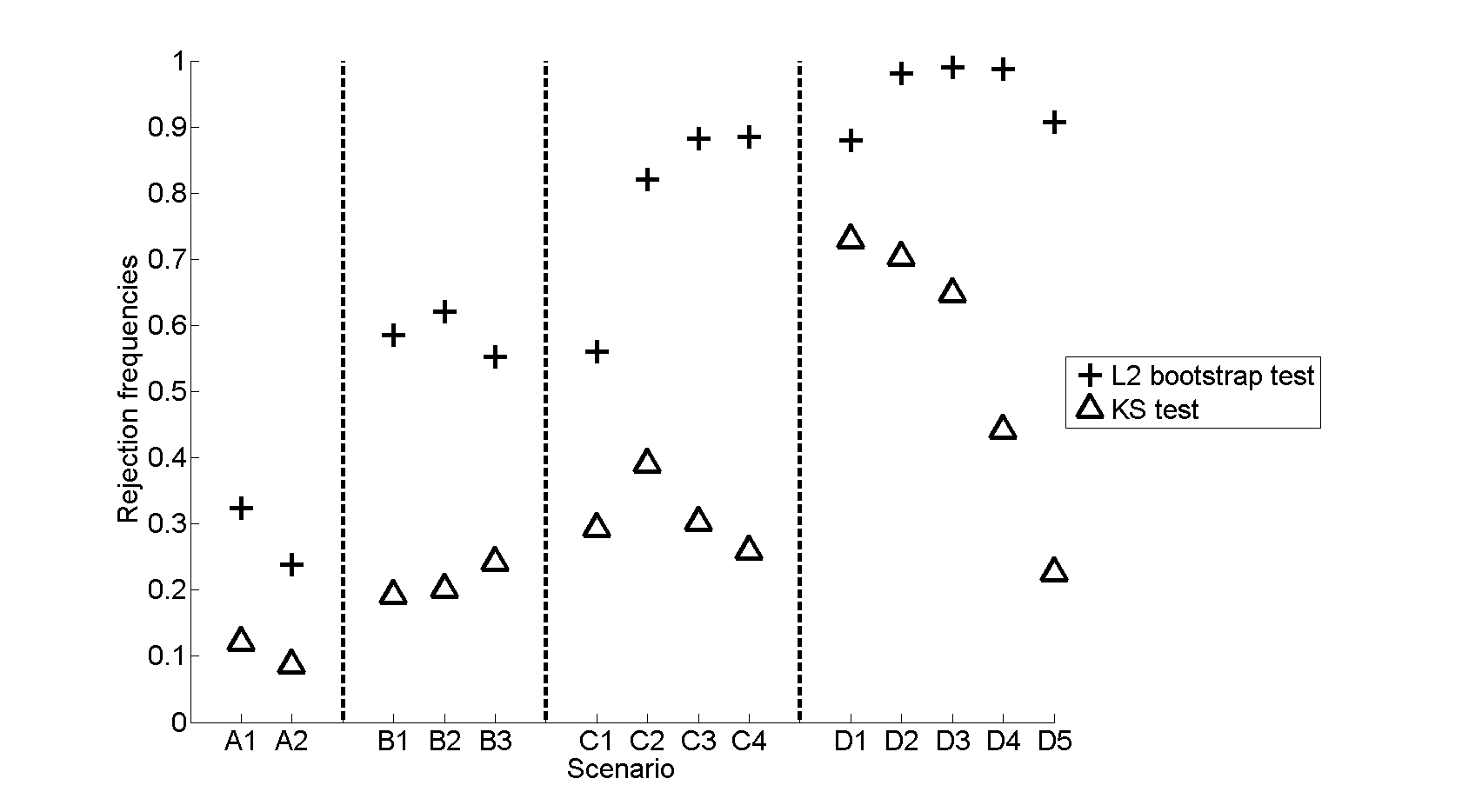}~~ \\ 
  \caption{ \label{alternative2} {\it Power of the test (\ref{test}) and the Kolmogorov-Smirnov test of \cite{prevet2012} at  5\% level for the model \eqref{alt2} under the scenarios A1-D5 from Table \ref{TableARd=0.1 mitboot}.} }
\end{center}
\end{figure}

\begin{figure}
\begin{center}
\includegraphics[width=1\textwidth,height=0.5\textwidth]{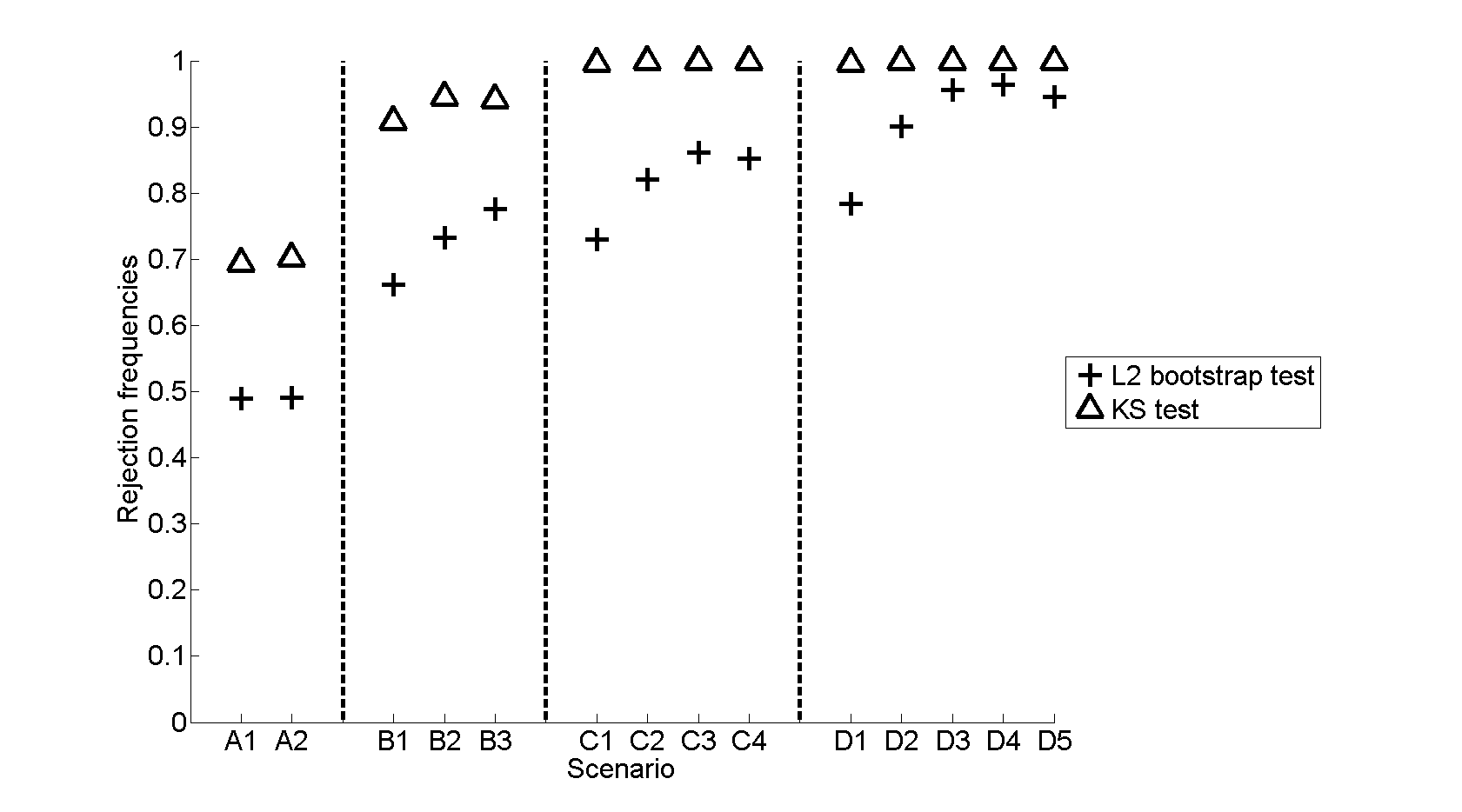}~~ \\
  \caption{ \label{alternative3} {\it Power of the test (\ref{test}) and the Kolmogorov-Smirnov test of \cite{prevet2012} at  5\% level for the model \eqref{alt3} under the scenarios A1-D5 from Table \ref{TableARd=0.1 mitboot}.} }
\end{center}
\end{figure}

In order to investigate the approximation of the nominal level we simulate data from the \\FARIMA($1,d,0$) model
\begin{eqnarray}
\label{FARIMA(1,d,0)}
(1-\phi B)(1-B)^d X_t=Z_t
\end{eqnarray}
and  the FARIMA($0,d,1$) process
\begin{eqnarray}
(1-B)^d X_t=(1+\theta B)Z_t
\label{FARIMA(0,d,1)}
\end{eqnarray}

for different values of $\phi,\theta$ and $d$ where the random variables $Z_t$ are independent standard normal distributed. The rejection probabilities  for the bootstrap test  \eqref{test} are displayed in Table \ref{TableARd=0.1 mitboot}--\ref{TableMAd=0.2 mitboot} where $d \in \{0.1,0.2\}$. We observe a very precise approximation of the nominal level in nearly all cases which is rather robust with respect to different choices of the parameter $M$ and $N$. 
\\
\\
In order to study the power of the test we consider the following alternatives
\begin{align}
(1-B)^{d}X_{t,T}&=Z_t+0.8\cos\left(1.5-\cos(4\pi t/T)\right)Z_{t-1}, \label{alt1} \\
\left(1-0.6\sin(4\pi t/T) B\right)(1-B)^{d}X_{t,T}&=Z_t, \label{alt2} \\
(1-B)^{d}X_{t,T}&=\sqrt{\sin(\pi t/T)}Z_t, \label{alt3}
\end{align}

where $d=0.2$. These kinds of alternatives were investigated by several authors in the context of locally stationary short-memory processes [see \cite{paparoditis2010} and \cite{dahlhaus1997}]. The rejection frequencies for the bootstrap test \eqref{test} are depicted in Figure \ref{alternative1}--\ref{alternative3} for different combinations of $T$ and $N$. Additionally, the results for the Kolmogorov-Smirnov approach of \cite{prevet2012} are presented. We observe that the new procedure clearly outperforms the test of \cite{prevet2012} for the models \eqref{alt1} and \eqref{alt2} while the Kolmogorov-Smirnov test works  better for the process \eqref{alt3}. In addition, we observe that the new decision rule is less sensitive with respect to different choices of $N$ than the test based on the Kolmogorov-Smirnov distance.

\bigskip

\begin{table}
\begin{center}
{\scriptsize
\begin{tabular}{|c|ccc|cc|cc|cc|cc|cc|}
\hline
&$$ & $$ & $$ &   \multicolumn{2}{c|}{$\phi=-0.9$}  &   \multicolumn{2}{c|}{$\phi=-0.5$}     & \multicolumn{2}{c|}{$\phi=0$} & \multicolumn{2}{c|}{$\phi=0.5$} &   \multicolumn{2}{c|}{$\phi=0.9$} \\
\hline
&$T$&$N$&$M$& $5 \%$ & $10 \%$   & $5 \%$ & $10 \%$   & $5 \%$ & $10 \%$  & $5 \%$ & $10\%$ & $5 \%$ & $10\%$ \\ \hline
A1&128&		16  &	8    &	.126&	.182	&	.072	&	.123	&	.036	&	.074	&	.073	&	.126	&	.084	&	.147\\
A2&128&		8    &	16  &	.140&	.200	&	.085	&	.132	&	.041	&	.090	&	.073	&	.128	&	.084	&	.118\\
B1&256&		32  &	8    &	.065&	.135	&	.064	&	.119	&	.042	&	.088	&	.062	&	.113	&	.075	&	.147\\
B2&256&		16  &	16  &	.080&	.132	&	.056	&	.108	&	.040	&	.082	&	.051	&	.109	&	.062	&	.109\\
B3&256&		8    &	32  &	.068&	.111	&	.045	&	.095	&	.046	&	.097	&	.072	&	.147	&	.049	&	.114\\
C1&512&		64  &	8    &	.054&	.109	&	.049	&	.106	&	.039	&	.089	&	.049	&	.114	&	.082	&	.134\\
C2&512&		32  &	16  &	.038&	.093	&	.043	&	.086	&	.039	&	.085	&	.059	&	.108	&	.065	&	.132\\
C3&512&		16  &	32  &	.061&	.095	&	.051	&	.102	&	.045	&	.081	&	.059	&	.109	&	.043	&	.104\\
C4&512&		8    &	64  &	.060&	.107	&	.053	&	.098	&	.045	&	.083	&	.060	&	.116	&	.042	&	.093\\
D1&1024&	128&	8    &	.039&	.104	&	.042	&	.093	&	.042	&	.085	&	.035	&	.093	&	.079	&	.132\\
D2&1024&	64  &	16  &	.053&	.104	&	.058	&	.097	&	.050	&	.110	&	.057	&	.101	&	.068	&	.126\\
D3&1024&	32  &	32  &	.033&	.076	&	.058	&	.114	&	.046	&	.086	&	.070	&	.107	&	.062	&	.115\\
D4&1024&	16  &	64  &	.046&	.089	&	.036	&	.091	&	.044	&	.084	&	.054	&	.109	&	.044	&	.099\\
D5&1024&	8    &	128&	.037&	.073	&	.041	&	.091	&	.041	&	.091	&	.061	&	.131	&	.045	&	.097\\
\hline
\end{tabular}}
\caption{\textit{Rejection probabilities of the bootstrap test \eqref{test} under $H_0$ for different choices of T,N and M. The data
was generated according to model \eqref{FARIMA(1,d,0)} with $d=0.1$ and different values for $\phi$.
}}
\label{TableARd=0.1 mitboot}
\end{center}
\begin{center}
{\scriptsize
\begin{tabular}{|c|ccc|cc|cc|cc|cc|cc|}
\hline
&$$ & $$ & $$ &   \multicolumn{2}{c|}{$\phi=-0.9$}  &   \multicolumn{2}{c|}{$\phi=-0.5$}     & \multicolumn{2}{c|}{$\phi=0$} & \multicolumn{2}{c|}{$\phi=0.5$} &   \multicolumn{2}{c|}{$\phi=0.9$} \\
\hline
&$T$&$N$&$M$& $5 \%$ & $10 \%$   & $5 \%$ & $10 \%$   & $5 \%$ & $10 \%$  & $5 \%$ & $10\%$ & $5 \%$ & $10\%$ \\ \hline
A1&128	&	16	&	8	&	.107	&	.164	&	.063	&	.114	&	.050	&	.108	&	.072	&	.121	&	.108	&	.166	\\
A2&128	&	8	&	16	&	.106	&	.160	&	.064	&	.118	&	.041	&	.085	&	.073	&	.124	&	.078	&	.138	\\
B1&256	&	32	&	8	&	.064	&	.123	&	.048	&	.104	&	.042	&	.094	&	.075	&	.131	&	.079	&	.137	\\
B2&256	&	16	&	16	&	.058	&	.125	&	.051	&	.101	&	.040	&	.101	&	.065	&	.112	&	.055	&	.116	\\
B3&256	&	8	&	32	&	.079	&	.124	&	.047	&	.089	&	.051	&	.091	&	.053	&	.106	&	.050	&	.105	\\
C1&512	&	64	&	8	&	.050	&	.093	&	.048	&	.090	&	.051	&	.103	&	.047	&	.104	&	.075	&	.133	\\
C2&512	&	32	&	16	&	.047	&	.104	&	.044	&	.087	&	.039	&	.085	&	.053	&	.109	&	.068	&	.124	\\
C3&512	&	16	&	32	&	.042	&	.097	&	.044	&	.087	&	.057	&	.106	&	.046	&	.105	&	.060	&	.104	\\
C4&512	&	8	&	64	&	.050	&	.102	&	.053	&	.101	&	.052	&	.088	&	.058	&	.121	&	.062	&	.114	\\
D1&1024&	128	&	8	&	.044	&	.090	&	.046	&	.102	&	.051	&	.107	&	.039	&	.092	&	.076	&	.140	\\
D2&1024&	64	&	16	&	.043	&	.082	&	.040	&	.088	&	.050	&	.098	&	.046	&	.098	&	.060	&	.106	\\
D3&1024&	32	&	32	&	.045	&	.089	&	.054	&	.097	&	.057	&	.103	&	.060	&	.104	&	.066	&	.115	\\
D4&1024&	16	&	64	&	.044	&	.087	&	.038	&	.087	&	.049	&	.094	&	.059	&	.106	&	.051	&	.101	\\
D5&1024&	8	&	128	&	.041	&	.082	&	.041	&	.089	&	.038	&	.086	&	.061	&	.103	&	.054	&	.103	\\
\hline
\end{tabular}
}
\caption{\textit{Rejection probabilities of the bootstrap test \eqref{test} under $H_0$ for different choices of T,N and M. The data
was generated according to model \eqref{FARIMA(1,d,0)} with $d=0.2$ and different values for $\phi$.
}}
\label{TableARd=0.2 mitboot}
\end{center}
\end{table}
\begin{table}%
\begin{center}
{\scriptsize
\begin{tabular}{|c|ccc|cc|cc|cc|cc|cc|}
\hline
&$$ & $$ & $$ &   \multicolumn{2}{c|}{$\theta=-0.9$}  &   \multicolumn{2}{c|}{$\theta=-0.5$}     & \multicolumn{2}{c|}{$\theta=0$} & \multicolumn{2}{c|}{$\theta=0.5$} &   \multicolumn{2}{c|}{$\theta=0.9$} \\
\hline
&$T$&$N$&$M$& $5 \%$ & $10 \%$   & $5 \%$ & $10 \%$   & $5 \%$ & $10 \%$  & $5 \%$ & $10\%$ & $5 \%$ & $10\%$ \\ \hline
A1&128	&	16	&	8	&	.072	&	.116	&	.054	&	.107	&	.041	&	.085	&	.044	&	.088	&	.077	&	.123	\\
A2&128	&	8	&	16	&	.068	&	.122	&	.054	&	.112	&	.059	&	.125	&	.073	&	.117	&	.070	&	.133	\\
B1&256	&	32	&	8	&	.045	&	.100	&	.060	&	.101	&	.041	&	.081	&	.042	&	.082	&	.036	&	.084	\\
B2&256	&	16	&	16	&	.053	&	.096	&	.058	&	.104	&	.045	&	.094	&	.045	&	.102	&	.060	&	.104	\\
B3&256	&	8	&	32	&	.064	&	.123	&	.057	&	.113	&	.049	&	.101	&	.042	&	.092	&	.061	&	.130	\\
C1&512	&	64	&	8	&	.043	&	.089	&	.043	&	.095	&	.044	&	.086	&	.045	&	.088	&	.041	&	.095	\\
C2&512	&	32	&	16	&	.046	&	.109	&	.067	&	.112	&	.052	&	.093	&	.051	&	.096	&	.043	&	.086	\\
C3&512	&	16	&	32	&	.048	&	.099	&	.055	&	.095	&	.062	&	.114	&	.050	&	.102	&	.051	&	.098	\\
C4&512	&	8	&	64	&	.038	&	.097	&	.055	&	.100	&	.047	&	.100	&	.046	&	.093	&	.042	&	.092	\\
D1&1024	&	128	&	8	&	.053	&	.103	&	.060	&	.099	&	.051	&	.099	&	.071	&	.118	&	.044	&	.094	\\
D2&1024	&	64	&	16	&	.044	&	.100	&	.062	&	.124	&	.048	&	.090	&	.068	&	.119	&	.042	&	.093	\\
D3&1024	&	32	&	32	&	.053	&	.107	&	.064	&	.116	&	.044	&	.082	&	.045	&	.094	&	.043	&	.098	\\
D4&1024	&	16	&	64	&	.044	&	.096	&	.038	&	.084	&	.042	&	.093	&	.045	&	.087	&	.042	&	.087	\\
D5&1024	&	8	&	128	&	.049	&	.109	&	.042	&	.085	&	.054	&	.109	&	.048	&	.083	&	.042	&	.096	\\
\hline
\end{tabular}}
\caption{\textit{Rejection probabilities of the bootstrap test \eqref{test} under $H_0$ for different choices of T,N and M. The data
was generated according to model \eqref{FARIMA(0,d,1)} with $d=0.1$ and different values for $\theta$.
}}
\label{TableMAd=0.1 mitboot}
\end{center}
\begin{center}
{\scriptsize
\begin{tabular}{|c|ccc|cc|cc|cc|cc|cc|}
\hline
&$$ & $$ & $$ &   \multicolumn{2}{c|}{$\theta=-0.9$}  &   \multicolumn{2}{c|}{$\theta=-0.5$}     & \multicolumn{2}{c|}{$\theta=0$} & \multicolumn{2}{c|}{$\theta=0.5$} &   \multicolumn{2}{c|}{$\theta=0.9$} \\
\hline
&$T$&$N$&$M$& $5 \%$ & $10 \%$   & $5 \%$ & $10 \%$   & $5 \%$ & $10 \%$  & $5 \%$ & $10\%$ & $5 \%$ & $10\%$ \\ \hline
A1&128	&	16	&	8	&	.068	&	.112	&	.060	&	.103	&	.030	&	.081	&	.060	&	.108	&	.053	&	.111	\\
A2&128	&	8	&	16	&	.060	&	.117	&	.051	&	.103	&	.061	&	.110	&	.062	&	.114	&	.068	&	.117	\\
B1&256	&	32	&	8	&	.059	&	.122	&	.048	&	.102	&	.045	&	.094	&	.038	&	.078	&	.040	&	.083	\\
B2&256	&	16	&	16	&	.053	&	.109	&	.041	&	.095	&	.047	&	.093	&	.040	&	.080	&	.048	&	.091	\\
B3&256	&	8	&	32	&	.060	&	.100	&	.048	&	.098	&	.057	&	.119	&	.050	&	.092	&	.061	&	.106	\\
C1&512	&	64	&	8	&	.059	&	.110	&	.064	&	.122	&	.052	&	.099	&	.056	&	.099	&	.055	&	.101	\\
C2&512	&	32	&	16	&	.060	&	.122	&	.044	&	.107	&	.041	&	.103	&	.043	&	.113	&	.046	&	.086	\\
C3&512	&	16	&	32	&	.061	&	.116	&	.056	&	.122	&	.049	&	.089	&	.046	&	.088	&	.052	&	.099	\\
C4&512	&	8	&	64	&	.056	&	.095	&	.057	&	.118	&	.057	&	.110	&	.047	&	.100	&	.055	&	.102	\\
D1&1024	&	128	&	8	&	.063	&	.125	&	.054	&	.102	&	.039	&	.086	&	.044	&	.101	&	.051	&	.098	\\
D2&1024	&	64	&	16	&	.051	&	.109	&	.061	&	.112	&	.047	&	.107	&	.056	&	.106	&	.047	&	.100	\\
D3&1024	&	32	&	32	&	.055	&	.092	&	.057	&	.111	&	.048	&	.095	&	.057	&	.106	&	.047	&	.119	\\
D4&1024	&	16	&	64	&	.065	&	.124	&	.061	&	.116	&	.043	&	.092	&	.048	&	.087	&	.049	&	.098	\\
D5&1024	&	8	&	128	&	.059	&	.116	&	.052	&	.095	&	.049	&	.093	&	.035	&	.075	&	.050	&	.115	\\
\hline
\end{tabular}
}
\caption{\textit{Rejection probabilities of the bootstrap test \eqref{test} under $H_0$ for different choices of T,N and M. The data
was generated according to model \eqref{FARIMA(0,d,1)} with $d=0.2$ and different values for $\theta$.
}}
\label{TableMAd=0.2 mitboot}
\end{center}
\end{table}


\phantom{.}


{\bf Acknowledgements.}
This work has been supported in part by the
Collaborative Research Center ``Statistical modeling of nonlinear
dynamic processes'' (SFB 823, Teilprojekt A1, C1) of the German Research Foundation (DFG).

\bibliographystyle{apalike}

\bibliography{presendet}

\section{Appendix: technical details}
\def\theequation{6.\arabic{equation}}
\setcounter{equation}{0}

In the following, we will state two results which will be central for the proof of the statements in Sections $3$ and $4$.
\begin{satz}
\label{thmapp1}
If Assumption \ref{ass1} is satisfied with $d_{\infty}<1/4$, the following statements are correct.
\begin{itemize}
\item[a)] 
\noindent
\vspace{-.7cm}
\bea
\E\big((\hat F_{1,T}, \hat F_{2,T})^T\big)=(F_{1,T},F_{2,T}+ d_{N,T})^T+\boldsymbol{\tilde{C}}_T +O\Big(\frac{1}{N^{1-4d_{\infty}}}\Big),
\eea
where the vector $\boldsymbol{\tilde{C}_T} \in \R^2$  is of order $O\big(N^2/T^2+\log(N)/(MN^{1-4d_{\infty}})\big)$. In particular, this term vanishes if the functions $\psi_l(u)$ are independent of $u$ for all $l \in \Z$.
\item[b)] 
\noindent
\vspace{-.7cm}
\bea
\Cov\big((\hat F_{1,T}, \hat F_{2,T})^T\big)=\boldsymbol{\tilde{\Sigma}_{T}}+O(T,d_{\infty}),
\eea
where 
\begin{eqnarray*}
\label{Otd}
O(T,d_{\infty}):=O\Big(\frac{\log (N)  }{ N^{1-8d_{\infty}} T}\Big) +O \Big( \frac{N^{2}}{T^3}   +\frac{N^{2+4d_{\infty}}}{T^3}1_{\{\frac{1}{8}  \leq d_{\infty} < \frac{1}{4} \}}\Big)
\end{eqnarray*}
and  $\tilde{\Sigma}_{T}$ is the same as the matrix $\Sigma$ in (\ref{3.3}) except that the integral $\int_{0}^{1}$ is replaced by $\frac{1}{TM}\sum_{j=1}^{M}$.  
\item[c)] If $d_{\infty}<1/8$ and $l_1,l_2 \in \N_0 $ with $l_1+l_2 \geq 3$, then we have
\bea
\cum(\sqrt{T}\hat F_{1,T} \boldsymbol{1}^T_{l_1}, \sqrt{T}\hat F_{2,T} \boldsymbol{1}^T_{l_2} )=O\big(T^{(1-\frac{l_1+l_2}{2})(1-8d_{\infty})}\big),
\eea
where $\boldsymbol{1}_{l_i} \in \R^{l_i}$ denotes a vector containing merely ones ($i=1,2$).
\end{itemize}
\end{satz}

It follows by the same arguments as given in Section $4$ of \cite{prevet2012}, that there exist  parameters $\hat \psi_{l,p}$ such that the bootstrap process $X_{t,T}^*$ defined in (\ref{bootproc}) can be represented as 
\be \label{bootstrapmainfty}
X_{t,T}^*=\sum_{l=0}^\infty \hat \psi_{l,p} Z_{t-l}^*,
\ee
where $Z_t^*$ are the innovations from part $2)$ of the bootstrap description. We now assume that the null hypothesis \eqref{H0} holds, and consider the process
\be
\label{xstern2}
X_{t,T,2}^*=\sum_{l=0}^\infty \psi_l Z_{t-l}^*,
\ee

where the coefficients $\psi_l=\psi_l(u)$ are the coefficients in (\ref{localappr}). We then define $\hat D_{T,2}^{2,*}$ as $\hat D_T^{2}$ in (\ref{2.8}) whereby the random variables $X_{t,T}$ are replaced by   $X_{t,T,2}^*$. The next theorem shows that the random variable $\hat D_T^{2,*}$ can be approximated by $\hat D_{T,2}^{2,*}$. 

\begin{satz}
\label{thmapp2}
Let $\alpha>0$ be fixed and denote with $A_T(\alpha)$ the event where $|\underline{\hat d} - \underline{d}| \leq \alpha/4$. If Assumption \ref{ass1} and the inequality
\be
\label{unglboots}
|\hat \psi_{l,p}-\psi_l| l^{1-\max(\underline{\hat d},\underline d)} \leq C \frac{p^4\log(T)^{3/2}}{\sqrt{T}} \quad \forall l \in \N
\ee
are satisfied, then
\bea
a) \quad && \E\Big((\hat D_T^{2,*}-\hat D_{T,2}^{2,*})1_{A_T(\alpha)}\Big)=O\Big(p^4\log(T)^{3/2}N^{4\underline{d}-1+\alpha}T^{-1/2}\Big),\\
b) \quad && \Var \Big((\hat D_T^{2,*}-\hat D_{T,2}^{2,*})1_{A_T(\alpha)}\Big)=O\Big({p^{8} \log(T)^{3}}{\log(N)^2N^{\max(8{\underline d}-1,0)+2\alpha}}T^{-2}\Big) .
\eea
\end{satz}


\subsection{Proof of Theorem \ref{thmapp1}:}
{\bf Proof of part a):} We define $\tilde t_j:=t_j-N/2+1$,
$ \tilde{\psi}_l(u_{j,p}):=\psi_l\Big(\frac{\tilde t_j+p}{T}\Big)$, $Z_{a,b}:=Z_{a-N/2+1+b}$ and write similar to \cite{supplementdetprevet2010}
\begin{eqnarray*}
\E[\hat{F}_{1,T}]&=&\E\Big( \frac{1}{T} \sum_{j=1}^{M} \sum_{k=1}^{\lfloor \frac{N}{2}\rfloor} I_N(u_j,\lambda_{k,N})^2  \Big)=\frac{1}{2}\E\Big( \frac{1}{T} \sum_{j=1}^{M} \sum_{k=-\lfloor \frac{N-1}{2}\rfloor }^{\lfloor \frac{N}{2}\rfloor} I_N(u_j,\lambda_{k,N})^2  \Big)\\
&=&\frac{1}{2T} \sum_{j=1}^M\sum_{k=-\lfloor \frac{N-1}{2}\rfloor }^{\lfloor \frac{N}{2}\rfloor} \frac{1}{(2 \pi N)^2} \sum_{p,q,r,s=0}^{N-1} \sum_{l,m,n,o=0}^{\infty}e^{-i(p-q+r-s)\lambda_{k,N}}\\
&&\psi_{\tilde t_j+p,T,l}\psi_{\tilde t_j+q,T,m}\psi_{\tilde t_j+r,T,n}  \psi_{\tilde t_j+s,T,o} \E[Z_{t_j,p-l}Z_{t_j,q-m}Z_{t_j,r-n}Z_{t_j,s-o}]\\
&=&E^1_{N,T}+E^2_{N,T}+A_{N,T}+B_{N,T},
\end{eqnarray*}
where we use the notation
\begin{eqnarray*}
E^1_{N,T}&=&\frac{1}{2T} \sum_{j=1}^M \sum_{k=-\lfloor \frac{N-1}{2}\rfloor }^{\lfloor \frac{N}{2}\rfloor} \frac{1}{(2 \pi N)^2} \sum_{p,q,r,s=0}^{N-1} \sum_{l,m,n,o=0}^{\infty}  \psi_l ( u_j)  \psi_m (u_j) \psi_n (u_j)  \psi_o (u_j) e^{-i(p-q+r-s)\lambda_{k,N}}  \\ 
&&\Big(\E[Z_{t_j,p-l}Z_{t_j,q-m}] \E[Z_{t_j,r-n}Z_{t_j,s-o}]+\E[Z_{t_j,p-l}Z_{t_j,s-o}] \E[Z_{t_j,q-m}Z_{t_j,r-n}] \Big),\\
E^2_{N,T}&=&
\frac{1}{2T} \sum_{j=1}^M \sum_{k=-\lfloor \frac{N-1}{2}\rfloor }^{\lfloor \frac{N}{2}\rfloor} \frac{1}{(2 \pi N)^2} \sum_{p,q,r,s=0}^{N-1}\sum_{l,m,n,o=0}^{\infty}\psi_l ( u_j)  \psi_m (u_j) \psi_n (u_j)  \psi_o (u_j)\\
&&e^{-i(p-q+r-s)\lambda_{k,N}} \E[Z_{t_j,p-l}Z_{t_j,r-n}] \E[Z_{t_j,q-m}Z_{t_j,s-o}],  \\
 A_{N,T}&=&\frac{1}{2T} \sum_{j=1}^M\sum_{k=-\lfloor \frac{N-1}{2}\rfloor }^{\lfloor \frac{N}{2}\rfloor} \frac{1}{(2 \pi N)^2 } \sum_{p,q,r,s=0}^{N-1} \sum_{l,m,n,o=0}^{\infty} 
e^{-i(p-q+r-s)\lambda_{k,N}}  \E[Z_{t_j,p-l}Z_{t_j,q-m}Z_{t_j,r-n}Z_{t_j,s-o}]\\
&& \Big \{ \Big(\tilde{\psi}_l(u_{j,p})- \psi_l(u_j) \Big)  \tilde{\psi}_m(u_{j,q})\tilde{\psi}_n(u_{j,r})\tilde{\psi}_o(u_{j,s}) +\psi_l(u_j) \Big(\tilde{\psi}_m(u_{j,q})- \psi_m(u_j) \Big)\tilde{\psi}_n(u_{j,r})\tilde{\psi}_o(u_{j,s})  \\
&&+\psi_l(u_j) \psi_m(u_j) \Big(\tilde{\psi}_n(u_{j,r})-\psi_n(u_j) \Big) \tilde{\psi}_o(u_{j,s})+\psi_l(u_j) \psi_m(u_j)\psi_n(u_j)  \Big(\tilde{\psi}_o(u_{j,s})- \psi_o(u_j)\Big)   \Big \} , \\
B_{N,T}&=&\frac{1}{2T} \sum_{j=1}^M\sum_{k=-\lfloor \frac{N-1}{2}\rfloor }^{\lfloor \frac{N}{2}\rfloor} \frac{1}{(2 \pi N)^2 } \sum_{p,q,r,s=0}^{N-1} \sum_{l,m,n,o=0}^{\infty} 
e^{-i(p-q+r-s)\lambda_{k,N}}  \E[Z_{t_j,p-l}Z_{t_j,q-m}Z_{t_j,r-n}Z_{t_j,s-o}]\\
&& \Big \{ \Big(\psi_{\tilde t_j+p,T,l}-\tilde{\psi}_l(u_{j,p})  \Big)  \tilde{\psi}_m(u_{j,q})\tilde{\psi}_n(u_{j,r})\tilde{\psi}_o(u_{j,s}) +\psi_{\tilde t_j+p,T,l} \Big(\psi_{\tilde t_j+q,T,m}-\tilde{\psi}_m(u_{j,q}) \Big)\tilde{\psi}_n(u_{j,r})\tilde{\psi}_o(u_{j,s})  \\
&&+\psi_{\tilde t_j+p,T,l}\psi_{\tilde t_j+q,T,m} \Big(\psi_{\tilde t_j+r,T,n}-\tilde{\psi}_n(u_{j,r}) \Big) \tilde{\psi}_o(u_{j,s})+\psi_{\tilde t_j+p,T,l}\psi_{\tilde t_j+q,T,m}\psi_{\tilde t_j+r,T,n}  \Big(\psi_{\tilde t_j+s,T,o}-\tilde{\psi}_o(u_{j,s}) \Big)   \Big \} .
\end{eqnarray*}

Note that $B_{N,T}$ corresponds to the error which occurs if the coefficients $\psi_{t,T,l}$ are replaced by $\psi_l(t/T)$ and that $A_{N,T}$ contains the approximation error  of $\psi_l(t/T)$ through $\psi_l(t_j/T)$ with $t_j$ denoting the midpoint of the $j$-th block. The  following four statements conclude the proof for $\E[\hat{F}_{1,T}]$
\begin{eqnarray}
E^1_{N,T}&=&\frac{1}{2\pi M}  \sum_{j=1}^M \int_{-\pi}^{\pi}  f^2(u_j,\lambda)\,d\lambda+ O \Big(\frac{1}{N^{1-4d_{\infty}}} \Big) ,
\label{E_{N,T}^1}\\
E^2_{N,T}&=&O \Big(\frac{1}{N^{1-4d_{\infty}}}\Big),
\label{E_{N,T}^2} \\
A_{N,T}&=&O\Big( \frac{\log N}{M N^{1-4d_{\infty}}} \Big)+ O\Big( \frac{N^2}{T^2}\Big),
\label{A_{N,T}} \\
B_{N,T}&=&O\Big(\frac{1}{T}\Big).
\label{B_{N,T}}
\end{eqnarray}
\textbf{Proof of \eqref{E_{N,T}^1}:} Without loss of generality, we only consider the first summand in $E^1_{N,T}$. Due to the independence of the random variables $Z_t$, we obtain that only those terms contribute to the sum where the conditions $0 \leq p=q+l-m \leq N-1$ and $0 \leq r=s+n-o \leq N-1$ are satisfied, which implies the inequality $\max\{|l-m|,|n-o|\} \leq N-1$. Thus, the first term in $E^1_{N,T}$ can  be expressed as
\begin{eqnarray*}
&&\frac{1}{2T} \sum_{j=1}^M \sum_{k=-\lfloor \frac{N-1}{2}\rfloor }^{\lfloor \frac{N}{2}\rfloor}   \frac{1}{(2\pi N)^2}  \sum_{\substack{l,m,n,o=0 \\  |l-m| \leq N-1 \\  |n-o| \leq N-1}}^{\infty} \psi_l(u_j) \psi_m(u_j) \psi_n(u_j) \psi_o(u_j) e^{-i(l-m+n-o)\lambda_{k,N}}\\
&&\hspace{.6cm}\times(N-|l-m|)(N-|n-o|) \\
&=&\frac{1}{2M} \sum_{j=1}^M \sum_{h=-1}^{1}  \frac{1}{(2\pi N)^2}  \sum_{\substack{l,m,n,o=0\\     |l-m| \leq N-1 \\  |n-o| \leq N-1 \\ l-m+n-o=hN}}^{\infty} \psi_l(u_j) \psi_m(u_j) \psi_n(u_j) \psi_{o}(u_j) (N-|l-m|)(N-|n-o|), 
\end{eqnarray*}
where we used the well known identity
\be \label{3.11}
 \sum_{k=0}^{N-1} \exp(-i \lambda_{k,N} r )= \begin{cases} N &  \mbox{if } r= Nh \mbox{ for some } h \in \Z \\ 0 &  \mbox{else, } \end{cases}
\ee
[note that we only have to consider three possible values of $h$ since $\max\{|l-m|,|n-o| \}\leq N-1$]. 
It is easy to see that 
\begin{eqnarray*}
E^1_{N,T}=E^1_{N,T,0}+E^1_{N,T,1}+E^1_{N,T,2}+E^1_{N,T,3}
\end{eqnarray*}
where
\begin{eqnarray}
E^1_{N,T,0}&=&\frac{1}{2M} \sum_{j=1}^M   \frac{1}{(2\pi)^2 }  \sum_{\substack{l,m,n,o=0\\   l-m+n-o=0}}^{\infty} \psi_l(u_j) \psi_m(u_j) \psi_n(u_j) \psi_{o}(u_j) =  \frac{1}{4\pi M}  \sum_{j=1}^M \int_{-\pi}^{\pi}  f^2(u_j,\lambda)\,d\lambda, 
\notag\\
E^1_{N,T,1}&=&\frac{1}{MN^2} \sum_{j=1}^M \sum_{h=-1}^{1}    \sum_{\substack{l,m,n,o=0\\     |l-m| \leq N-1 \\  |n-o| \leq N-1 \\ l-m+n-o=hN}}^{\infty} \psi_l(u_j) \psi_m(u_j) \psi_n(u_j) \psi_{o}(u_j) 
\label{E_{N,T}^1 1.Fehler}\\
&&\hspace{1cm}\times (-N|l-m|-N|n-o|+|l-m||n-o|), 
\notag \\
E^1_{N,T,2}&=&\frac{1}{M} \sum_{j=1}^M \sum_{h \in \{-1,1 \}}^{}    \sum_{\substack{l,m,n,o=0\\     |l-m| \leq N-1 \\  |n-o| \leq N-1 \\ l-m+n-o=hN}}^{\infty} \psi_l(u_j) \psi_m(u_j) \psi_n(u_j) \psi_{o}(u_j) ,
\label{E_{N,T}^1 2.Fehler}\\
E^1_{N,T,3}&=&\frac{1}{M} \sum_{j=1}^M    \sum_{\substack{l,m,n,o=0\\     N \leq |l-m| \\ N \leq  |n-o|  \\ l-m+n-o=0}}^{\infty} \psi_l(u_j) \psi_m(u_j) \psi_n(u_j) \psi_{o}(u_j). 
\label{E_{N,T}^1 3.Fehler}
\end{eqnarray}
In order to complete the proof of ($\ref{E_{N,T}^1}$), it therefore suffices to demonstrate that the last three expressions are of order $O \big(\frac{1}{N^{1-4d_{\infty}}} \big)$. We commence with (\ref{E_{N,T}^1 1.Fehler}). 
Setting $o=l-m+n-hN \geq 0$ and using  (\ref{apprpsi}), it follows that there exists a constant $C\in \R$ such that 
\begin{eqnarray}
E^{1}_{N,T,1} &\leq & \frac{C}{MN} \sum_{j=1}^M \sum_{h =-1}^{1}   \sum_{\substack {l,m,n=1 \\1\leq|l-m-hN| \leq N-1 \\ 1 \leq l-m+n-hN}}^{\infty}\frac{1}{l^{1-d_{\infty}}}\frac{1}{m^{1-d_{\infty}}}\frac{1}{n^{1-d_{\infty}}}\frac{|l-m-hN|}{(l-m+n-hN)^{1-d_{\infty}}} 
\label{69a}
\end{eqnarray}
(note that all terms where one of the variables $l,m,n$ or $l-m+n-hN$ vanishes are of smaller or the same order). This argument will be employed continuously throughout this proof without mentioning it explicitly. Note that the summand $|l-m|$ does not occur in the numerator of the above expression due to the symmetry of $|l-m|$ and $|n-o|$ in \eqref{E_{N,T}^1 1.Fehler}, and that $C \in \R^{+}$ denotes a universal constant throughout the whole proof.  Setting $z:=l-m-hN$, we obtain $|z|=|l-m-hN| \leq N-1$ and the expression on the right hand side of (\ref{69a}) can be written as
\begin{eqnarray*}
  && \frac{C}{N}  \sum_{h =-1}^{1} \sum_{\substack{z \in \Z \\ 1 \leq |z| \leq N-1}}^{}  \sum_{\substack {m,n=1  \\ 1 \leq z+m+hN \\ 1 \leq n+z}}^{\infty}\frac{1}{(z+m+hN)^{1-d_{\infty}}}\frac{1}{m^{1-d_{\infty}}}\frac{1}{n^{1-d_{\infty}}}\frac{|z|}{(n+z)^{1-d_{\infty}}} \\
  &\stackrel{(\ref{eq1})}{\lesssim}& \frac{1}{N}  \sum_{h =-1}^{1}   \sum_{\substack{z \in \Z \\ 1 \leq |z| \leq N-1}}^{}\frac{|z|^{2d_{\infty}}}{|z+hN|^{1-2d_{\infty}}} 
 \lesssim  \frac{1}{N^{1-2d_{\infty}}} \sum_{h =-1}^{1}   \sum_{\substack{z \in \Z \\ 1 \leq |z| \leq N-1}}\frac{1}{|z+hN|^{1-2d_{\infty}}} \lesssim \frac{1}{N^{1-4d_{\infty}}} ,
\end{eqnarray*}
where $a_n \lesssim b_n$ means that $a_n/b_n$ is bounded by some finite constant for all $n \in \N$. By using \eqref{apprpsi},  \eqref{3.11}, \eqref{eq1} and similar arguments,  we obtain that (\ref{E_{N,T}^1 2.Fehler}) is bounded by
\begin{eqnarray*}
&&  \sum_{h \in \{-1,1 \}}^{}   \sum_{\substack{z \in \Z \\ 1 \leq |z| \leq N-1}}^{}\frac{1}{|z+hN|^{1-2d_{\infty}}} \frac{1}{|z|^{1-2d_{\infty}}}\stackrel{(\ref{Gleichung 3.196(3)}),(\ref{eq1})}{\lesssim}  \frac{1}{N^{1-4d_{\infty}}},
\end{eqnarray*}
and since (\ref{E_{N,T}^1 3.Fehler}) is shown analogously, we therefore conclude the proof of \eqref{E_{N,T}^1}. \\

{\bf Proof of (\ref{E_{N,T}^2})}: The result follows by similar arguments as used in the treatment of (\ref{E_{N,T}^1 1.Fehler})--(\ref{E_{N,T}^1 3.Fehler}).
\\
\\
{\bf Proof of (\ref{A_{N,T}}):} Without loss of generality, we only consider the first summand and replace $\tilde{\psi}_m(u_{j,q})\tilde{\psi}_n(u_{j,r})\tilde{\psi}_o(u_{j,s})$ by $\psi_m(u_j)\psi_n(u_j)\psi_o(u_j)$
[the error due to this replacement is negligible, which follows by analogous arguments as given for the term $A_{N,T,1}^{(2)}$ at a later stage of this proof]. Due to the independence of the random variables $Z_t$, we obtain the sum of  three terms [compare the definition of $E_{N,T}^1$ for the first two summands and the definition of $E_{N,T}^2$ for the third one] and we restrict ourselves to the first one, i.e. we only consider
\bea
A_{N,T,1}&:=&\frac{1}{2T} \sum_{j=1}^M\sum_{k=-\lfloor \frac{N-1}{2}\rfloor }^{\lfloor \frac{N}{2}\rfloor} \frac{1}{(2 \pi N)^2 } \sum_{p,q,r,s=0}^{N-1} \sum_{l,m,n,o=0}^{\infty} 
e^{-i(p-q+r-s)\lambda_{k,N}} \\  &&\E[Z_{t_j,p-l}Z_{t_j,q-m}]\E[Z_{t_j,r-n}Z_{t_j,s-o}]
 \Big(\tilde{\psi}_l(u_{j,p})- \psi_l(u_j) \Big)  \psi_m(u_j)\psi_n(u_j)\psi_o(u_j) .
\eea
Using a Taylor expansion, we can write
\begin{eqnarray*}
 \tilde{\psi}_l(u_{j,p})- \psi_l(u_j) =\psi_{l}^{'}(u_j)\Big(\frac{-N/2+1+p}{T} \Big)+ \frac{\psi^{''}(\eta_{l,j,p})}{2}\Big(\frac{-N/2+1+p}{T} \Big)^2
\end{eqnarray*}
with $\eta_{l,j,p} \in (u_j-N/(2T),u_j+N/(2T))$, and therefore $A_{N,T,1}$ splits into two terms which  will be denoted by $A_{N,T,1}^{(1)}$ and $A_{N,T,1}^{(2)}$ in the following discussion. We start with the treatment of the first summand. Employing the independence of the innovations we obtain that the indices corresponding to non vanishing terms must satisfy $q=p+m-l$ and $n=o+r-s$. Applying  (\ref{3.11}) afterwards  yields $0 \leq m=l+r-s - hN \mbox{ with }h \in \{-1,0,1 \}$ and this combined with \eqref{apprpsi} and (\ref{2.1b}) implies
\begin{eqnarray*}
A_{N,T,1}^{(1)}&\lesssim &\frac{1}{ N^2}\sum_{h=-1}^{1}\sum_{r,s=0}^{N-1} \sum_{\substack{l,o=1 \\1 \leq  o+r-s \\ 1 \leq l+r-s-hN \\ 0 \leq |r-s-hN| \leq N-1}}^{\infty}
\frac{\log (l)}{l^{1-d_{\infty}}}\frac{1}{(l+r-s-hN)^{1-d_{\infty}}}\frac{1}{(o+r-s)^{1-d_{\infty}}}\frac{1}{o^{1-d_{\infty}}} \\
&&\bigg|\sum_{\substack{p=0 \\0 \leq  p+r-s-hN \leq N-1}}^{N-1}\Big(\frac{-N/2+1+p}{T}\Big)\bigg|. 
\end{eqnarray*}
We restrict ourselves to the cases $|r-s| \geq 1$ and $|r-s-hN| \geq 1$ since the remaining terms are of smaller order. A straightforward calculation yields
\begin{eqnarray}
\label{Nq}
&&\bigg|\sum_{\substack{p=0 \\ 0 \leq p+q \leq N-1 }}^{N-1}\Big(\frac{-N/2+1+p}{T} \Big)\bigg| \\
&=& \frac{N}{2T} \times 1_{\{q=0\}}+\min\left(\frac{ N|q| }{T},\frac{(N-|q|)|q|}{T}\right) O(1) \times 1_{\{1 \leq |q| \leq N-1\}}
\end{eqnarray}
and by using the second summand it follows that $A_{N,T,1}^{(1)}$ is bounded by
\begin{eqnarray*}
&&\frac{C}{NT}\sum_{h=-1}^{1} \sum_{r,s=0}^{N-1} \sum_{\substack{l,o=1 \\1 \leq  o+r-s \\ 1 \leq l+r-s-hN \\ 1 \leq |r-s-hN| \leq N-1}}^{\infty}
\frac{\log (l)}{l^{1-d_{\infty}}}\frac{1}{(l+r-s-hN)^{1-d_{\infty}}}\frac{1}{(o+r-s)^{1-d_{\infty}}}\frac{1}{o^{1-d_{\infty}}}|r-s-hN|  \\
&\lesssim&\frac{\log (N)}{N^{1-2d_{\infty}}T} \sum_{r,s=0}^{N-1} \sum_{\substack{o=1 \\1 \leq  o+r-s }}^{\infty}\frac{1}{(o+r-s)^{1-d_{\infty}}}\frac{1}{o^{1-d_{\infty}}}  \lesssim \frac{\log (N)}{N^{1-2d_{\infty}}T} \sum_{r,s=0}^{N-1}\frac{1}{|r-s|^{1-2d_{\infty}}}\lesssim  \frac{\log (N)}{MN^{1-4d_{\infty}}}
\end{eqnarray*}
where we used Lemma \ref{Lemma 7.1}c) and \ref{Lemma 7.1}b) for the first and second inequality, respectively. Next, we show that $A_{N,T,1}^{(2)}$ is of order $O (N^2T^{-2})$, and for this reason we choose $\epsilon>0$ such that $1-4d_{\infty}-\epsilon >0$. Using \eqref{apprpsi}, (\ref{2.1b}), (\ref{2.1c}) and $h \in \{-1,0,1 \}$,  the claim then follows by a further application of Lemma \ref{Lemma 7.1}b)
\begin{eqnarray*}
&&\frac{C}{ N^2} \sum_{r,s=0}^{N-1} \sum_{\substack{l,o=1 \\1 \leq  o+r-s \\ 1 \leq l+r-s-hN }}^{\infty}
\frac{\log^2 (l)}{l^{1-d_{\infty}}}\frac{1}{(l+r-s-hN)^{1-d_{\infty}}}\frac{1}{(o+r-s)^{1-d_{\infty}}}\frac{1}{o^{1-d_{\infty}}}\\
&&\times\sum_{\substack{p=0 \\0 \leq  p+r-s-hN \leq N-1}}^{N-1}\Big(\frac{-N/2+1+p}{T}\Big)^2 \\
&\lesssim &\frac{N}{ T^2}\sum_{r,s=0}^{N-1} \sum_{\substack{l,o=1 \\1 \leq  o+r-s \\ 1 \leq l+r-s-hN }}^{\infty}
\frac{1}{l^{1-d_{\infty}-\epsilon}}\frac{1}{(l+r-s-hN)^{1-d_{\infty}}}\frac{1}{(o+r-s)^{1-d_{\infty}}}\frac{1}{o^{1-d_{\infty}}} \\
&\lesssim&\frac{N}{ T^2}\sum_{r,s=0}^{N-1} \frac{1}{|r-s-hN|^{1-2d_{\infty}-\epsilon}}\frac{1}{|r-s|^{1-2d_{\infty}}}\stackrel{}{\lesssim} \frac{N^2}{ T^2}.
\end{eqnarray*}

\textbf{Proof of (\ref{B_{N,T}}):} The statement follows from  (\ref{apprbed}) and similar arguments as given in the proofs of (\ref{E_{N,T}^1}) and (\ref{A_{N,T}}). 
\\
\\
In order to proof the assertion for $\hat{F}_{2,T}$, one proceeds in the same way and the details are omitted for the sake of brevity. However it turns out that the expression corresponding to $E^2_{N,T}$ does not vanish in this case and there appears an additional  bias, which will be denoted by $d_{N,T}$.
\\
\\
{\bf Proof of part b)} We restrict ourselves to the proof of  
\begin{eqnarray*}
 &&\Var ( \hat F_{1,T} )=\frac{1}{4}\Var \Bigl(\frac{1}{T} \sum_{j=1}^{M}\sum_{k=-\lfloor \frac{N-1}{2}\rfloor }^{\lfloor \frac{N}{2}\rfloor} I_N(u_{j}, \lambda_{k,N})^2 \Bigr) =
\frac{5}{\pi}\frac{1}{TM} \sum_{j=1}^M   \int_{-\pi}^{\pi}f^4(u_{j}, \lambda)\,d\lambda+ O(T,d_{\infty})
\end{eqnarray*}
and recall the definition of the remainder
\begin{eqnarray*}
O(T,d_{\infty})=O\Big(\frac{\log (N)  }{ N^{1-8d_{\infty}} T}\Big) +O \Big( \frac{N^{2}}{T^3}   +\frac{N^{2+4d_{\infty}}}{T^3}1_{\{\frac{1}{8}  \leq d_{\infty} < \frac{1}{4} \}}\Big).
\end{eqnarray*}
All other statements can be verified completely analogously and the details are omitted. By combining the arguments from the proof of part a) and from
\cite{supplementdetprevet2010}, we obtain that
\bea
\Var \Bigl(\frac{1}{T} \sum_{j=1}^{M}\sum_{k=-\lfloor \frac{N-1}{2}\rfloor }^{\lfloor \frac{N}{2}\rfloor} I_N(u_{j}, \lambda_{k,N})^2 \Bigr)=32V^{*}(\nu^{'})+8V^{*}(\nu^{''})+ O(T,d_{\infty}),
\eea
where
\begin{eqnarray*}
V^{*}(\nu^{'})&=&\frac{1}{T^2} \sum_{j_1,j_2=1}^M  \sum_{k_1,k_2=-\lfloor \frac{N-1}{2}\rfloor }^{\lfloor \frac{N}{2}\rfloor} \frac{1}{(2\pi N)^4}\sum_{p_1,q_1,r_1,s_1=0}^{N-1} \sum_{p_2,q_2,r_2,s_2=0}^{N-1} \sum_{v_1,w_1,x_1,y_1=0}^{\infty} \sum_{v_2,w_2,x_2,y_2 =0}^{\infty}  \\ 
&& \psi_{v_1}(u_{j_1})\psi_{w_1}(u_{j_1})\psi_{x_1}(u_{j_1})\psi_{y_1}(u_{j_1})\psi_{v_2}(u_{j_2})\psi_{w_2}(u_{j_2}) \psi_{x_2}(u_{j_2})\psi_{y_2}(u_{j_2}) \\
&&e^{-i (p_1-q_1+r_1-s_1)  \lambda_{k_1}}e^{-i(p_2-q_2+r_2-s_2)\lambda_{k_2}} \\
&& \E[ Z_{t_{j_1},p_1-v_1} Z_{t_{j_1},q_1-w_1}] \E[ Z_{t_{j_1},r_1-x_1} Z_{t_{j_2},p_2-v_2}] \E [Z_{t_{j_1},s_1-y_1}Z_{t_{j_2},q_2-w_2}] \E [Z_{t_{j_2},r_2-x_2}Z_{t_{j_2},s_2-y_2}]
\end{eqnarray*}
and
\begin{eqnarray*}
V^{*}(\nu^{''})
&=&\frac{1}{T^2} \sum_{j_1,j_2=1}^M  \sum_{k_1,k_2=-\lfloor \frac{N-1}{2}\rfloor }^{\lfloor \frac{N}{2}\rfloor} \frac{1}{(2\pi N)^4} \sum_{p_1,q_1,r_1,s_1=0}^{N-1} \sum_{p_2,q_2,r_2,s_2=0}^{N-1}\sum_{v_1,w_1,x_1,y_1=0}^{\infty} \sum_{v_2,w_2,x_2,y_2 =0}^{\infty}  \\ 
&& \psi_{v_1}(u_{j_1})\psi_{w_1}(u_{j_1})\psi_{x_1}(u_{j_1})\psi_{y_1}(u_{j_1})\psi_{v_2}(u_{j_2})\psi_{w_2}(u_{j_2}) \psi_{x_2}(u_{j_2})\psi_{y_2}(u_{j_2}) \\
&&e^{-i (p_1-q_1+r_1-s_1)  \lambda_{k_1}}e^{-i(p_2-q_2+r_2-s_2)\lambda_{k_2}}   \E[ Z_{t_{j_1},p_1-v_1}  Z_{t_{j_2},p_2-v_2}] \E[ Z_{t_{j_1},r_1-x_1} Z_{t_{j_2},r_2-x_2}] \\
&&\E [Z_{t_{j_1},q_1-w_1}Z_{t_{j_2},q_2-w_2}] \E [Z_{t_{j_1},s_1-y_1}Z_{t_{j_2},s_2-y_2}  ] .
\end{eqnarray*}

We start with $V^{*}(\nu^{'})$. Because of the independence of the random variables $Z_t$, the restrictions
$p_1 = q_1+v_1-w_1$, $p_2= r_1+ v_2-x_1+({j_1}-{j_2})N$, $q_2 =s_1+w_2-y_1+(j_1-j_2)N$ and $s_2 = r_2  +y_2-x_2$ are necessary for a non vanishing term. Consider $h_1,h_2 \in \{ -1,0,1 \}$ and sum over $k_1, k_2$ by using $(\ref{3.11})$. Then, $V^{*}(\nu^{'})$ can be written as
\begin{eqnarray*}
V^{*}(\nu^{'})&=&\frac{1}{T^2} \sum_{j_1,j_2=1}^M  \sum_{k_1,k_2=-\lfloor \frac{N-1}{2}\rfloor }^{\lfloor \frac{N}{2}\rfloor} \frac{1}{(2\pi N)^4}
\sum_{q_1,r_1,s_1,r_2=0}^{N-1}
\sum_{\substack{v_1,w_1,x_1,y_1=0 \\  0 \leq  q_1+v_1-w_1  \leq N-1 }}^{\infty} \sum_{\substack{v_2,w_2,x_2,y_2 =0 \\ 0 \leq r_1+ v_2-x_1+({j_1}-{j_2})N \leq N-1 \\ 0 \leq  s_1+w_2-y_1+(j_1-j_2)N \leq N-1\\   0 \leq r_2  +y_2-x_2 \leq N-1}}^{\infty} \\ 
&& \psi_{v_1}(u_{j_1})\psi_{w_1}(u_{j_1})\psi_{x_1}(u_{j_1})\psi_{y_1}(u_{j_1})\psi_{v_2}(u_{j_2})\psi_{w_2}(u_{j_2}) \psi_{x_2}(u_{j_2})\psi_{y_2}(u_{j_2}) \\
&&e^{-i (v_1-w_1+r_1-s_1)  \lambda_{k_1}}e^{-i( r_1+ v_2-x_1-s_1-w_2+y_1 -y_2+x_2)\lambda_{k_2}}\\
&=&\frac{1}{M^2} \sum_{j_1,j_2=1}^M \sum_{h_1,h_2=-1}^{1}\frac{1}{(2\pi N)^4}\sum_{\substack{r_1,s_1=0  }}^{N-1} \sum_{\substack{v_1,w_1,x_1,y_1=0  \\ v_1-w_1+r_1-s_1=h_1N \\ 0 \leq |v_1-w_1| \leq N-1 \\ }}^{\infty} \sum_{\substack{v_2,w_2,x_2,y_2 =0\\ r_1+ v_2-x_1-s_1-w_2+y_1 -y_2+x_2 =h_2N\\ 0 \leq r_1+ v_2-x_1+({j_1}-{j_2})N \leq N-1\\ 0 \leq  s_1+w_2-y_1+(j_1-j_2)N \leq N-1 \\ 0 \leq |y_2-x_2| \leq N-1}}^{\infty} \\ 
&& \psi_{v_1}(u_{j_1})\psi_{w_1}(u_{j_1})\psi_{x_1}(u_{j_1})\psi_{y_1}(u_{j_1})\psi_{v_2}(u_{j_2})\psi_{w_2}(u_{j_2}) \psi_{x_2}(u_{j_2})\psi_{y_2}(u_{j_2}) \\
&&(N^2-N|y_2-x_2|-N|v_1-w_1|+|v_1-w_1|y_2-x_2| ).
\end{eqnarray*}
An application of \eqref{apprpsi} yields similar to the proof of part a) that the above expression is of order $O \big( \frac{1}{N^{1-8d_{\infty}}T}\big)$, if
\begin{itemize}
\item[(i)]  $h_1,h_2 \in \{ -1,1 \}$ [compare \eqref{E_{N,T}^1 2.Fehler}],
\item[(ii)]  $j_1 \not= j_2$ [we prove this claim in Lemma \ref{j_1 ungleich j_2} in the appendix since this kind of restriction did not occur in the proof of part a)],
\item[(iii)] we drop $-N|y_2-x_2|-N|v_1-w_1|+|v_1-w_1|y_2-x_2|$  [compare \eqref{E_{N,T}^1 1.Fehler}],
\item[(iv)]  we drop $0 \leq |v_1-w_1| \leq N-1 $ and $ 0 \leq |y_2-x_2| \leq N-1$ [compare \eqref{E_{N,T}^1 3.Fehler}],
\item[(v)] we drop
$0 \leq r_1+ v_2-x_1+({j_1}-{j_2})N \leq N-1$ and $0 \leq  s_1+w_2-y_1+(j_1-j_2)N \leq N-1$ [compare \eqref{E_{N,T}^1 1.Fehler}].
\end{itemize}
By rearranging the equation $v_1-w_1+r_1-s_1=0$ to $0 \leq s_1 =r_1+ v_1-w_1 \leq N-1$, it follows
\begin{eqnarray*}
V^{*}(\nu^{'}) &=&\frac{1}{M^2N^2} \sum_{j_1=1}^M   \frac{1}{(2\pi )^4} \sum_{\substack{r_1=0 }}^{N-1}\sum_{\substack{v_1,w_1,x_1,y_1=0 \\
0 \leq r_1+v_1-w_1 \leq N-1 }}^{\infty} \sum_{\substack{v_2,w_2,x_2,y_2 =0 \\ w_1-v_1+ v_2-x_1-w_2+y_1 -y_2+x_2 =0}}^{\infty} \\ 
&& \psi_{v_1}(u_{j_1})\psi_{w_1}(u_{j_1})\psi_{x_1}(u_{j_1})\psi_{y_1}(u_{j_1})\psi_{v_2}(u_{j_1})\psi_{w_2}(u_{j_1}) \psi_{x_2}(u_{j_1})\psi_{y_2}(u_{j_1}) + O(T,d_{\infty})\\
&=&\frac{1}{M^2N} \sum_{j_1=1}^M  \frac{1}{(2\pi )^4}\sum_{\substack{v_1,w_1,x_1,y_1=0 }}^{\infty} \sum_{\substack{v_2,w_2,x_2,y_2 =0 \\ w_1-v_1+ v_2-w_2+y_1-x_1 +x_2-y_2=0}}^{\infty}  \\ 
&& \psi_{v_1}(u_{j_1})\psi_{w_1}(u_{j_1})\psi_{x_1}(u_{j_1})\psi_{y_1}(u_{j_1})\psi_{v_2}(u_{j_1})\psi_{w_2}(u_{j_1}) \psi_{x_2}(u_{j_1})\psi_{y_2}(u_{j_1})  + O(T,d_{\infty}) \\
&=& \frac{1}{TM}\frac{1}{2\pi} \sum_{j_1=1}^M \int_{-\pi}^{\pi}f^4(u_{j_1}, \lambda) \,d\lambda+ O(T,d_{\infty}).
\end{eqnarray*}
By using the same techniques as in $V^{*}(\nu^{'})$, we obtain  
\begin{eqnarray*}
V^{*}(\nu^{''})
&=&\frac{1}{M^2} \sum_{j_1=1}^M \sum_{h_1=-1}^{1}\frac{1}{(2\pi N)^4} \sum_{\substack{q_2,r_2,s_2=0 \\ 0 \leq  q_2-r_2+s_2+h_1N \leq N-1}}^{N-1} \sum_{v_1,w_1,x_1,y_1=0}^{\infty} \sum_{\substack{v_2,w_2,x_2,y_2 =0 \\ v_1-v_2+w_2-w_1+x_1-x_2+y_2-y_1 =0}}^{\infty}  \\ 
&& \psi_{v_1}(u_{j_1})\psi_{w_1}(u_{j_1})\psi_{x_1}(u_{j_1})\psi_{y_1}(u_{j_1})\psi_{v_2}(u_{j_1})\psi_{w_2}(u_{j_1}) \psi_{x_2}(u_{j_1})\psi_{y_2}(u_{j_1}) +O(T,d_{\infty}).
\end{eqnarray*}
Note that, in contrast to the term $V^{*}(\nu^{'})$, the cases where $h_1 \in \{-1,1\}$ do not 
vanish. In fact, using the three equalities 
\begin{eqnarray*}
&&\sum_{\substack{q_2,r_2,s_2=0 \\ 0 \leq  q_2-r_2+s_2 \leq N-1}}^{N-1}= \frac{2}{3}N^3+O( N^2) ,
\sum_{\substack{q_2,r_2,s_2=0 \\ 0 \leq  q_2-r_2+s_2+N \leq N-1}}^{N-1}= \frac{1}{6}N^3+O( N^2), \\
&&\sum_{\substack{q_2,r_2,s_2=0 \\ 0 \leq  q_2-r_2+s_2-N \leq N-1}}^{N-1}= \frac{1}{6}N^3+O( N^2) 
\end{eqnarray*}
we deduce
\begin{eqnarray}
V^{*}(\nu^{''})&=&\Big(\frac{2}{3}+\frac{1}{6}+\frac{1}{6}\Big)\frac{1}{M^2N} \sum_{j_1=1}^M   \frac{1}{(2\pi )^4} \sum_{v_1,w_1,x_1,y_1=0}^{\infty} \sum_{\substack{v_2,w_2,x_2,y_2 =0 \\ v_1-v_2+w_2-w_1+x_1-x_2+y_2-y_1 =0}}^{\infty}  \nonumber \\ 
&& \psi_{v_1}(u_{j_1})\psi_{w_1}(u_{j_1})\psi_{x_1}(u_{j_1})\psi_{y_1}(u_{j_1})\psi_{v_2}(u_{j_1})\psi_{w_2}(u_{j_1}) \psi_{x_2}(u_{j_1})\psi_{y_2}(u_{j_1})   +O(T,d_{\infty}) \nonumber  \\
&=&\frac{1}{TM} \frac{1}{2\pi}\sum_{j_1=1}^M \int_{-\pi}^{\pi}f^4(u_{j_1}, \lambda) \,d\lambda+O(T,d_{\infty}). \label{varianzunterschied}
\end{eqnarray}
{\bf Proof of part c)} Exemplarily we consider the case $l_2=0$  and $l:=l_1 \geq 3 $.
The other cases can be treated similarly with an additional amount of notation.
Following the same lines as in the proof of Theorem 3.1 in \cite{detprevet2010}, it is sufficient to choose an arbitrary indecomposable partition 
\bea
\Big \{(Z_{t_{i_1},a_{1}-v_{1}}   Z_{t_{i_2},a_{2}-w_{1}}), (Z_{t_{i_3},a_{3}-x_{1}}  Z_{t_{i_4},a_{4}-y_{1}}),...,(Z_{t_{i_{4l-1}},a_{4l-1}-x_{l}}  Z_{t_{i_{4l}},a_{4l}-y_{l}} ) \Big \}
\eea
of the table
\be
\begin{matrix} \label{table}
Z_{t_{j_1},p_{1}-g_{1}} & Z_{t_{j_1},q_{1}-m_{1}} & Z_{t_{j_1},r_{1}-n_{1}} & Z_{t_{j_1},s_{1}-o_{1}} \\
\vdots & \vdots & \vdots & \vdots \\
Z_{t_{j_l},p_{l}-g_{l}} & Z_{t_{j_l},q_{l}-m_{l}} & Z_{t_{j_l},r_{l}-n_{l}} & Z_{t_{j_l},s_{l}-o_{l}} \\
\end{matrix}
\ee

[see \cite{brillinger1981}] and to treat the term
\begin{eqnarray*}
&&\frac{T^{\frac{l}{2} }  }{T^l} \frac{1} {N^{2l}}\sum_{j_1,...,j_l=1}^{M} \sum_{k_1, \ldots ,k_l=-\lfloor \frac{N-1}{2}\rfloor }^{\lfloor \frac{N}{2}\rfloor}  \sum_{p_1,...,s_{l}=0}^{N-1}\sum_{v_{1},w_{1},x_{1},y_{1}=0}^{\infty} \cdots \sum_{v_{l},w_{l},x_{l},y_{l}=0}^{\infty} \\
&& \tilde{\psi}_{v_{1}}(u_{j_1,p_1})\tilde{\psi}_{w_{1}}(u_{j_1,q_1})  \tilde{\psi}_{x_{1}}(u_{j_1,r_1})\tilde{\psi}_{y_{1}}(u_{j_1,s_1}) \tilde{\psi}_{v_{2}}(u_{j_2,p_2})\tilde{\psi}_{w_{2}}(u_{j_2,q_2})  \tilde{\psi}_{x_{2}}(u_{j_2,r_2})\tilde{\psi}_{y_{2}}(u_{j_2,s_2})\\
&&  \cdots \\
&& \tilde{\psi}_{v_{l}}(u_{j_l,p_l})\tilde{\psi}_{w_{l}}(u_{j_l,q_l})  \tilde{\psi}_{x_{l}}(u_{j_l,r_l})\tilde{\psi}_{y_{l}}(u_{j_l,s_l})  e^{-i (p_1-q_1+r_1-s_1) \lambda_{k_1}} \cdots e^{-i(p_l-q_l+r_l-s_l) \lambda_{k_l} } \\
&&\E[Z_{t_{i_1},a_{1}-v_{1}}   Z_{t_{i_2},a_{2}-w_{1}}  ]  \E[ Z_{t_{i_3},a_{3}-x_{1}}  Z_{t_{i_4},a_{4}-y_{1}}]  \E[Z_{t_{i_5},a_{5}-v_{2}}   Z_{t_{i_6},a_{6}-w_{2}} ]\E[Z_{t_{i_7},a_{7}-x_{2}}  Z_{t_{i_8},a_{8}-y_{2}}] \\
&& \hdots \\
&& \E[Z_{t_{i_{4l-3}},a_{4l-3}-v_{l}}   Z_{t_{i_{4l-2}},a_{4l-2}-w_{l}} ]\E[Z_{t_{i_{4l-1}},a_{4l-1}-x_{l}}  Z_{t_{i_{4l}},a_{4l}-y_{l}} ] 
\end{eqnarray*}
with $\{a_{1},a_2 \ldots, a_{4l} \} \in \{p_1, \ldots, p_l, q_1, \ldots, q_l, r_1, \ldots, r_l, s_1, \ldots, s_l \}$, $a_i \not = a_j$ for $i\not = j$,  
$\{ {i_1}, {i_2} ,\ldots ,{i_{4l}}\} \in  \{  {j_1},{j_2}, \ldots   {j_{l}}\}$, and $|\{ {i_1}, {i_2} ,\ldots ,{i_{4l}}\}|=l$. We now discuss the conditions which yield a contribution different from $0$ in this sum. Note that some of the $i_k$ are equal to each other and we will therefore write $j_1,...,j_l$ for the $l$ different values and consider $i_k$ as a function depending on $j_1,...,j_l$. Using the independence of the random variables $Z_t$ and summing with respect to $k_1, \ldots, k_l$, the conditions
\begin{eqnarray}
&&a_{4m+1}-a_{4m+2}+w_{m+1}-v_{m+1}+(i_{4m+1}-i_{4m+2})N=0 \mbox{ for } m=0,\ldots, l-1 ,
\label{611}\\
&&a_{4m+3}-a_{4m+4}+y_{m+1}-x_{m+1}+(i_{4m+3}-i_{4m+4})N=0 \mbox{ for } m=0,\ldots, l-1 ,
\label{612} \\
&&p_i-q_i+r_i-s_i=h_iN \mbox{ for } i =1,2, \ldots ,l \mbox{ and } h_i \in \{ -1,0,1\}
\label{Bedingung *}
\end{eqnarray}
follow. Rearranging the equations in (\ref{Bedingung *}) for a variable and plugging them into the $l$ equations \eqref{611} (where in every equation only \underline{one} variable is replaced) yields, due to the indecomposability of the partition and $v_{m+1}, x_{m+1} \geq 0 $, that the conditions
\begin{eqnarray}
\label{2l Gleichungen}
\begin{array}{lllll}
(1) &0 \leq v_1= \tilde{a}_{1}-\tilde{a}_{2}+\tilde{a}_{3}-\tilde{a}_4+w_{1}+(i_1-i_2+h_1)N ,    \\
(2) &0 \leq v_2=\tilde{a}_{7}-\tilde{a}_{8}+\tilde{a}_9-\tilde{a}_{10}+w_{2}+(i_5-i_6+h_2)N ,  \\
&\vdots \\
(l) &0 \leq v_l= \tilde{a}_{6l-5}-\tilde{a}_{6l-4}+\tilde{a}_{6l-3}-\tilde{a}_{6l-2}+w_{l}+(i_{4l-3}-i_{4l-2}+h_l)N , \\
(l+1) &0 \leq x_1= \tilde{a}_{5}-\tilde{a}_{6}+y_{1}+(i_{3}-i_{4})N, \\
(l+2)&0 \leq  x_{2}= \tilde{a}_{11}-\tilde{a}_{12}+y_{2}+(i_{7}-i_{8})N , \\
&\vdots \\
 (2l) &0 \leq  x_{l}= \tilde{a}_{6l-1}-\tilde{a}_{6l}+y_{l}+(i_{4l-1}-i_{4l})N
\end{array}
\end{eqnarray}
must hold, where $\{\tilde{a}_1, \tilde{a}_2, \ldots, \tilde{a}_{6l} \} \in \{p_1,  \ldots, s_l \}$ and $ |\{\tilde{a}_1, \tilde{a}_2, \ldots, \tilde{a}_{6l} \}| = 3l$. By employing  \eqref{apprpsi}, we can bound the above expression up to a constant by
\begin{eqnarray*}
& &\frac{1}{M^l}\frac{T^{\frac{l}{2}}  }{ N^{2l}}\sum_{j_1,...,j_l=1}^{M}\sum_{h_1, \ldots, h_l=-1}^{1}\sum_{\tilde{a}_1,...,\tilde{a}_{6l}=0}^{N-1}\\
&&\sum_{\substack{w_{1},y_{1}=1 \\ \tilde{a}_{1}-\tilde{a}_{2}+\tilde{a}_{3}-\tilde{a}_4+w_{1}+(i_1-i_2+h_1)N \geq 1 \\ \tilde{a}_{5}-\tilde{a}_{6}+y_{1}+(i_{3}-i_{4})N \geq 1}}^{\infty} \cdots \sum_{\substack{w_{l},y_{l}=1 \\ \tilde{a}_{6l-5}-\tilde{a}_{6l-4}+\tilde{a}_{6l-3}-\tilde{a}_{6l-2}+w_{l}+(i_{4l-3}-i_{4l-2}+h_l)N \geq 1 \\ \tilde{a}_{6l-1}-\tilde{a}_{6l}+y_{l}+(i_{4l-1}-i_{4l})N \geq 1 }}^{\infty} \\
&&\frac{1}{(\tilde{a}_{1}-\tilde{a}_{2}+\tilde{a}_{3}-\tilde{a}_4+w_{1}+(i_1-i_2+h_1)N)^{1-d_{\infty}}}\frac{1}{w_{1}^{1-d_{\infty}}}\frac{1}{(\tilde{a}_{5}-\tilde{a}_{6}+y_{1}+(i_{3}-i_{4})N
)^{1-d_{\infty}}}\frac{1}{y_{1}^{1-d_{\infty}}}\\
&& \frac{1}{(\tilde{a}_{7}-\tilde{a}_{8}+\tilde{a}_9-\tilde{a}_{10}+w_{2}+(i_5-i_6+h_2)N  )^{1-d_{\infty}}}\frac{1}{w_{2}^{1-d_{\infty}}}\frac{1}{(\tilde{a}_{11}-\tilde{a}_{12}+y_{2}+(i_{7}-i_{8})N)^{1-d_{\infty}}}\frac{1}{y_{2}^{1-d_{\infty}}}\\
&&  \cdots \\
&&\frac{1}{(\tilde{a}_{6l-5}-\tilde{a}_{6l-4}+\tilde{a}_{6l-3}-\tilde{a}_{6l-2}+w_{l}+(i_{4l-3}-i_{4l-2}+h_l)N)^{1-d_{\infty}}}\frac{1}{w_{l}^{1-d_{\infty}}}\\
&&\frac{1}{(\tilde{a}_{6l-1}-\tilde{a}_{6l}+y_{l}+(i_{4l-1}-i_{4l})N)^{1-d_{\infty}}}\frac{1}{y_{l}^{1-d_{\infty}}}.
\end{eqnarray*}

Using Lemma \ref{Lemma 7.1}b) in the appendix, this term can be (up to a constant) bounded by
\begin{eqnarray*}
&&\frac{1}{M^l}\frac{T^{\frac{l}{2}}  }{ N^{2l}}\sum_{j_1,...,j_l=1}^{M}\sum_{h_1, \ldots, h_l=-1}^{1}\sum_{\substack{\tilde{a}_1, \tilde{a}_2, \ldots, \tilde{a}_6=0 \\ |\tilde{a}_{1}-\tilde{a}_{2}+\tilde{a}_{3}-\tilde{a}_4+(i_1-i_2+h_1)N| \geq 1 \\ |\tilde{a}_{5}-\tilde{a}_{6}+(i_{3}-i_{4})N| \geq 1}}^{N-1} \cdots \sum_{\substack{\tilde{a}_{6l-5}, \tilde{a}_{6l-4}, \ldots, \tilde{a}_{6l}=0 \\ |\tilde{a}_{6l-5}-\tilde{a}_{6l-4}+\tilde{a}_{6l-3}-\tilde{a}_{6l-2}+(i_{4l-3}-i_{4l-2}+h_l)N| \geq 1 \\ |\tilde{a}_{6l-1}-\tilde{a}_{6l}+(i_{4l-1}-i_{4l})N| \geq 1 }}^{N-1} \\
&&\frac{1}{|\tilde{a}_{1}-\tilde{a}_{2}+\tilde{a}_{3}-\tilde{a}_4+(i_1-i_2+h_1)N|^{1-2d_{\infty}}}\frac{1}{|\tilde{a}_{5}-\tilde{a}_{6}+(i_{3}-i_{4})N
|^{1-2d_{\infty}}}\\
&& \frac{1}{|\tilde{a}_{7}-\tilde{a}_{8}+\tilde{a}_9-\tilde{a}_{10}+(i_5-i_6+h_2)N|^{1-2d_{\infty}}}\frac{1}{|\tilde{a}_{11}-\tilde{a}_{12}+(i_{7}-i_{8})N|^{1-2d_{\infty}}}\\
&&  \cdots \\
&&\frac{1}{|\tilde{a}_{6l-5}-\tilde{a}_{6l-4}+\tilde{a}_{6l-3}-\tilde{a}_{6l-2}+(i_{4l-3}-i_{4l-2}+h_l)N|^{1-2d_{\infty}}}\frac{1}{|\tilde{a}_{6l-1}-\tilde{a}_{6l}+(i_{4l-1}-i_{4l})N|^{1-2d_{\infty}}}.
\end{eqnarray*}

We now assume without loss of generality that
\begin{eqnarray}
 & &\tilde{a}_{6m+1}-\tilde{a}_{6m+2}+\tilde{a}_{6m+3}-\tilde{a}_{6m+4}+(i_{4m+1}-i_{4m+2}+h_{m+1})N \geq 1  ,
 \notag \\
& & \tilde{a}_{6m+5}-\tilde{a}_{6m+6}+(i_{4m+3}-i_{4m+4})N  \geq 1 
\label{Bedingungen Kumulante}
\end{eqnarray}
holds for $m=0,1,2,\ldots, l-1$ (the more general case follows analogously with an additional amount of notation). In this case the absolute values in the above expression can be skipped. It follows, as in  \cite{detprevet2010}, that the conditions on the $\tilde a_i$ imply that, if $i_1$ is chosen, there are only finitely many possible choices for $i_k$, $k=2,...,l$. Thus it suffices to consider  the following sum
\begin{eqnarray*}
&&\frac{1}{M^{l-1}}\frac{T^{\frac{l}{2}}  }{ N^{2l}}\sum_{\substack{\tilde{a}_1, \tilde{a}_2, \ldots, \tilde{a}_6=0 \\ \tilde{a}_{1}-\tilde{a}_{2}+\tilde{a}_{3}-\tilde{a}_4+C_1N \geq 1 \\ \tilde{a}_{5}-\tilde{a}_{6}+C_{l+1}N \geq 1}}^{N-1} \cdots \sum_{\substack{\tilde{a}_{6l-5}, \tilde{a}_{6l-4}, \ldots, \tilde{a}_{6l}=0 \\ \tilde{a}_{6l-5}-\tilde{a}_{6l-4}+\tilde{a}_{6l-3}-\tilde{a}_{6l-2}+C_lN \geq 1 \\ \tilde{a}_{6l-1}-\tilde{a}_{6l}+C_{2l}N \geq 1 }}^{N-1} \\
&&\frac{1}{(\tilde{a}_{1}-\tilde{a}_{2}+\tilde{a}_{3}-\tilde{a}_4+C_1N)^{1-2d_{\infty}}}
\frac{1}{(\tilde{a}_{5}-\tilde{a}_{6}+C_{l+1}N)^{1-2d_{\infty}}}\\
&& \frac{1}{(\tilde{a}_{7}-\tilde{a}_{8}+\tilde{a}_9-\tilde{a}_{10}+C_2N  )^{1-2d_{\infty}}}
\frac{1}{(\tilde{a}_{11}-\tilde{a}_{12}+C_{l+2}N)^{1-2d_{\infty}}}\\
&&  \cdots \\
&&\frac{1}{(\tilde{a}_{6l-5}-\tilde{a}_{6l-4}+\tilde{a}_{6l-3}-\tilde{a}_{6l-2}+C_lN)^{1-2d_{\infty}}}
\frac{1}{(\tilde{a}_{6l-1}-\tilde{a}_{6l}+C_{2l}N)^{1-2d_{\infty}}}
\end{eqnarray*}
with $C_1, C_2, \ldots, C_l \in \{ -1,0 \ldots, M \}$ and $C_{l+1}, C_{l+2}, \ldots, C_{2l} \in \{0,1 \ldots, M-1\}$ (because of  (\ref{Bedingungen Kumulante}) and  $\tilde{a}_i \in \{0,1,2, \ldots, N-1 \}$ there are no other possible values for $C_i$). 
We remind that (due to the indecomposability of the partition) the $2l$-fractions inside the addend are hooked. This means that for two different fractions there exists a chain of fractions (starting with the first considered fraction and ending with the second one), such that in every element of the chain there exists at least one element $\tilde a_i$ which also occurs in the consecutive fraction.
 We will perform a summation in a particular way and in order to illustrate this, we consider the first two fractions and assume that  $\tilde{a}_1$ and $\tilde{a}_6$ are (up to a the algebraic sign) the same. We distinguish two cases.
\begin{itemize}
\item[(i)]  If $\tilde{a}_1 =\tilde{a}_6$, we obtain from Lemma \ref{Lemma 7.1}a) that
\begin{eqnarray}
&&\sum_{\substack{\tilde{a}_1=0 \\ \tilde{a}_{1}-\tilde{a}_{2}+\tilde{a}_{3}-\tilde{a}_4+C_1N \geq 1 \\ \tilde{a}_{5}-\tilde{a}_{6} +C_{l+1}N\geq 1}}^{N-1}\frac{1}{(\tilde{a}_{1}-\tilde{a}_{2}+\tilde{a}_{3}-\tilde{a}_4+C_1N)^{1-2d_{\infty}}}\frac{1}{(\tilde{a}_{5}-\tilde{a}_{6}+C_{l+1}N)^{1-2d_{\infty}}} \nonumber \\
&& \lesssim \frac{1}{(-\tilde{a}_{2}+\tilde{a}_{3}-\tilde{a}_4+\tilde{a}_5+(C_1+C_{l+1})N)^{1-4d_{\infty}}} 
\label{Fall 1 ohne 2d}\\
&  & \lesssim \frac{T^{2d_{\infty}}}{(-\tilde{a}_{2}+\tilde{a}_{3}-\tilde{a}_4+\tilde{a}_5+(C_1+C_{l+1})N)^{1-2d_{\infty}}}.
\label{Fall 1 mit 2d}
\end{eqnarray}
Furthermore we have $-\tilde{a}_{2}+\tilde{a}_{3}-\tilde{a}_4+\tilde{a}_5+(C_1+C_{l+1})N \geq 2$ which follows from the conditions  $\tilde{a}_{1}-\tilde{a}_{2}+\tilde{a}_{3}-\tilde{a}_4+C_1N \geq 1$  and $\tilde{a}_{5}-\tilde{a}_{6}+C_{l+1}N= \tilde{a}_5-\tilde{a}_1+C_{l+1}N \geq 1$.  
\item[(ii)] If $\tilde{a}_1=-\tilde{a}_6$ and  $-\tilde{a}_{2}+\tilde{a}_{3}-\tilde{a}_4-\tilde{a}_5+(C_1-C_{l+1})N \not =0$, it follows from Lemma \ref{Lemma 7.1}b) that
\begin{eqnarray}
&&\sum_{\substack{\tilde{a}_1=0 \\ \tilde{a}_{1}-\tilde{a}_{2}+\tilde{a}_{3}-\tilde{a}_4+C_1N \geq 1 \\ \tilde{a}_{5}-\tilde{a}_{6}+C_{l+1}N \geq 1}}^{N-1}\frac{1}{(\tilde{a}_{1}-\tilde{a}_{2}+\tilde{a}_{3}-\tilde{a}_4+C_1N)^{1-2d_{\infty}}}\frac{1}{(\tilde{a}_{5}-\tilde{a}_{6}+C_{l+1}N)^{1-2d_{\infty}}} \nonumber \\
&  & \lesssim \frac{1}{|-\tilde{a}_{2}+\tilde{a}_{3}-\tilde{a}_4-\tilde{a}_5+(C_1-C_{l+1})N|^{1-4d_{\infty}}}
\label{Fall 2 ohne 2d} \\
&  & \lesssim \frac{T^{2d_{\infty}}}{|-\tilde{a}_{2}+\tilde{a}_{3}-\tilde{a}_4-\tilde{a}_5+(C_1-C_{l+1})N|^{1-2d_{\infty}}}.
\label{Fall 2 mit 2d}
\end{eqnarray}
\end{itemize}
In both cases, it is possible that variables cancel out, for example if $\tilde{a}_4=\tilde{a}_5$  and $\tilde{a}_3=\tilde{a}_5$ in the first and second case, respectively.
 We apply (\ref{Fall 1 ohne 2d})--(\ref{Fall 2 mit 2d}) in total $2l-2$-times. 
In the first  $2l-4$-applications, we use $(\ref{Fall 1 mit 2d})$ and $(\ref{Fall 2 mit 2d})$ (depending on the algebraic sign of the variable which appears in both fractions) and in the  $(2l-3)$th and $(2l-2)$th application we employ $(\ref{Fall 1 ohne 2d})$ and $(\ref{Fall 2 ohne 2d})$.
We furthermore assume that  $h$ variables cancel out while utilizing these inequalities. Afterwards $3l-(2l-2)-h= l+2-h$ variables remain with $0 \leq h \leq l$, namely $\tilde a_{6l-1}$, $\tilde a_{6l}$ and $l-h$ other variables with values in $\{ 0,1,2, \ldots, N-1\}$. Denoting these $l-h$ variables with  $b_1, b_2, \ldots, b_{l-h}$  we obtain
\begin{eqnarray*}
&&\frac{1}{M^{l-1}}\frac{T^{\frac{l}{2}}  }{ N^{2l}}\sum_{\substack{\tilde{a}_1, \tilde{a}_2, \ldots, \tilde{a}_6=0 \\ \tilde{a}_{1}-\tilde{a}_{2}+\tilde{a}_{3}-\tilde{a}_4+C_1N \geq 1 \\ \tilde{a}_{5}-\tilde{a}_{6}+C_{l+1}N \geq 1}}^{N-1} \cdots \sum_{\substack{\tilde{a}_{6l-5}, \tilde{a}_{6l-4}, \ldots, \tilde{a}_{6l}=0 \\ \tilde{a}_{6l-5}-\tilde{a}_{6l-4}+\tilde{a}_{6l-3}-\tilde{a}_{6l-2}+C_lN \geq 1 \\ \tilde{a}_{6l-1}-\tilde{a}_{6l}+C_{2l}N \geq 1 }}^{N-1} \\
&&\frac{1}{(\tilde{a}_{1}-\tilde{a}_{2}+\tilde{a}_{3}-\tilde{a}_4+C_1N)^{1-2d_{\infty}}}
\frac{1}{(\tilde{a}_{5}-\tilde{a}_{6}+C_{l+1}N)^{1-2d_{\infty}}} \\
&&\frac{1}{(\tilde{a}_{7}-\tilde{a}_{8}+\tilde{a}_9-\tilde{a}_{10}+C_2N  )^{1-2d_{\infty}}}
\frac{1}{(\tilde{a}_{11}-\tilde{a}_{12}+C_{l+2}N)^{1-2d_{\infty}}}\\
&&  \cdots \\
&&\frac{1}{(\tilde{a}_{6l-5}-\tilde{a}_{6l-4}+\tilde{a}_{6l-3}-\tilde{a}_{6l-2}+C_lN)^{1-2d_{\infty}}}
\frac{1}{(\tilde{a}_{6l-1}-\tilde{a}_{6l}+C_{2l}N)^{1-2d_{\infty}}} \\
& \lesssim& \frac{1}{M^{l-1}}\frac{T^{\frac{l}{2}}  }{ N^{2l}}\sum_{\substack{\tilde{a}_{6l-1}, \tilde{a}_{6l}=0\\ \tilde{a}_{6l-1}-\tilde{a}_{6l}+C_{2l}N \geq 1 }}^{N-1}\sum_{\substack{ b_1, b_2, \ldots, b_{l-h} =0 \\1 \leq  |\tilde{a}_{6l-1}-\tilde{a}_{6l}+\sum_{j=1}^{l-h} (-1)^{k_j}2b_{j} +\sum_{j=1}^{2l-1}(-1)^{k_j} C_jN|}}^{N-1}\\
&&\frac{N^hT^{(2l-4)2d_{\infty}}}{ |\tilde{a}_{6l-1}-\tilde{a}_{6l}+\sum_{j=1}^{l-h} (-1)^{k_j}2b_{j} +\sum_{j=1}^{2l-1}(-1)^{k_j} C_jN|^{1-6d_{\infty}}}  \frac{1}{(\tilde{a}_{6l-1}-\tilde{a}_{6l}+C_{2l}N)^{1-2d_{\infty}}} 
\end{eqnarray*}
with some $k_j \in \{0,1 \}$. We first consider the case $h=l$. If $\sum_{j=1}^{2l-1}(-1)^{k_j} C_jN=0$ and $C_{2l}=0$, it follows that the above term equals
\begin{eqnarray*}
&& \frac{1}{M^{l-1}}\frac{T^{\frac{l}{2}} N^lT^{(2l-4)2d_{\infty}} }{ N^{2l}}\sum_{\substack{\tilde{a}_{6l-1},\tilde{a}_{6l}=0 \\1 \leq \tilde{a}_{6l-1}-\tilde{a}_{6l} } }^{N-1}  \frac{1}{(\tilde{a}_{6l-1}-\tilde{a}_{6l})^{2-8d_{\infty}}}\lesssim  \frac{1}{M^{l-1}}\frac{T^{\frac{l}{2}} N^{l+1}T^{(2l-4)2d_{\infty}} }{ N^{2l}} =   T^{(1-\frac{l}{2})(1-8d_{\infty})}. 
\end{eqnarray*}
If $ |\sum_{j=1}^{2l-1}(-1)^{k_j} C_lN | \geq 1$ or $C_{2l} =1 $, we apply Lemma \ref{Lemma 7.1}a) and b) in order to obtain the same upper bound (it can be shown that, in this case, there appears an  additional factor $N^{1-8d_{\infty}}$   in the denominator, so the corresponding term is, in fact, of smaller order). The same upper bound holds for  $h \leq l-1$. 
  $\hfill \Box$
\subsection{Proof of Remark \ref{riemannstattintegral}: }
If we replace $\hat F_{1,T}$ by the corresponding integrated version $\tilde{F}_{1,T}=\frac{1}{4\pi M }\sum_{j=1}^{M}\int_{-\pi}^{\pi} I_N(u_j,\lambda)^2 \,d\lambda$, the derivation of the asymptotic variance  can be carried out almost analogously as in the proof of Theorem \ref{thmapp1}b) except that the term, where the variable $h_1$ in $V^*(v^{''})$ equals $-1$ or $1$, does not occur, because for the integrated version one can use
\be \label{3.11integral}
 \frac{1}{2\pi}\int_{-\pi}^\pi \exp(-i \lambda r )d\lambda= \begin{cases} 1 &  \mbox{if } r= 0, \\ 0 &  \mbox{else, } \end{cases}
\ee

(for $r \in \Z$) instead of \eqref{3.11}. Therefore, in the integrated case, we obtain
\begin{eqnarray*}
V^{*}(\nu^{''})&=&\frac{2}{3}\frac{1}{M^2N} \sum_{j_1=1}^M \frac{1}{(4\pi)^2}  \frac{1}{(2\pi )^2} \sum_{v_1,w_1,x_1,y_1=0}^{\infty} \sum_{\substack{v_2,w_2,x_2,y_2 =0 \\ v_1-v_2+w_2-w_1+x_1-x_2+y_2-y_1 =0}}^{\infty}   \\ 
&& \psi_{v_1}(u_{j_1})\psi_{w_1}(u_{j_1})\psi_{x_1}(u_{j_1})\psi_{y_1}(u_{j_1})\psi_{v_2}(u_{j_1})\psi_{w_2}(u_{j_1}) \psi_{x_2}(u_{j_1})\psi_{y_2}(u_{j_1})   +O(T,d_{\infty})   \\
&=&\frac{2}{3}\frac{1}{TM} \frac{1}{8\pi}\sum_{j_1=1}^M \int_{-\pi}^{\pi}f^4(u_{j_1}, \lambda) \,d\lambda+O(T,d_{\infty})
\end{eqnarray*}
instead of \eqref{varianzunterschied} and we recall that the order $O(T,d_{\infty})$ is defined in (\ref{Otd}). This yields that the asymptotic variance of $\sqrt{T} \tilde{F}_{1,T} $ is $\frac{14}{3\pi}\int_{-\pi}^\pi\int_0^1 f^4(u,\lambda) du d\lambda$ and does not coincide with the asymptotic variance of $\sqrt{T} \hat{F}_{1,T} $.

\subsection{Proof of Theorem \ref{thmapp2}:}

{\bf Proof of part a):} We define $\hat F_{1,T}^{*}$ and $\hat F_{1,T,2}^{*}$   as $\hat F_{1,T}$ where the observed data $X_{t,T}$ are replaced by $X_{t,T}^*$ and $X_{t,T,2}^*$, respectively. By using \eqref{bootstrapmainfty} and writing $I_N^*(u,\lambda)$ for the bootstrap analogue of $ I_N(u,\lambda)$, we then get
\bea
\E((\hat F_{1,T}^{*}-\hat F_{1,T,2}^{*})1_{A_T(\alpha)} | X_{1,T},...,X_{T,T})
&=& \frac{1}{2T} \sum_{j=1}^M\sum_{k=-\lfloor \frac{N-1}{2}\rfloor }^{\lfloor \frac{N}{2}\rfloor} \frac{1}{(2 \pi N)^2} \sum_{p,q,r,s=0}^{N-1} \sum_{l,m,n,o=0}^{\infty}e^{-i(p-q+r-s)\lambda_{k,N}}\\
&&\hat \psi_{l,m,n,o,p}1_{A_T(\alpha)}
 \E[Z^{*}_{t_j,p-l}Z^{*}_{t_j,q-m}Z^{*}_{t_j,r-n}Z^{*}_{t_j,s-o}]
\eea

[compare the first set of equalities in the proof of Theorem \ref{thmapp1} a)], where $\hat \psi_{l,m,n,o,p}=\hat \psi_{l,p}\hat \psi_{m,p}\hat \psi_{n,p}\hat \psi_{o,p}-\psi_l\psi_m\psi_n\psi_o$. By using the decomposition
\bea
\hat \psi_{l,m,n,o,p}&=&(\hat \psi_{l,p}-\psi_l)\hat \psi_{m,p}\hat \psi_{n,p}\hat \psi_{o,p}+\psi_l(\hat \psi_{m,p}-\psi_m)\hat \psi_{n,p}\hat \psi_{o,p} \\
&&+\psi_l\psi_m(\hat \psi_{n,p}-\psi_n)\hat \psi_{o,p} +\psi_l\psi_m\psi_n(\hat \psi_{o,p}-\psi_o)
\eea

the above expression splits into four terms and for the sake of brevity we only consider the first one. The other cases are treated similarly. As in the proof of Theorem \ref{thmapp1} a) we then obtain terms $E_{N,T}^{1,*}$ and $E_{N,T}^{2,*}$ which are defined as $E_{N,T}^{1}$, $E_{N,T}^{2}$ where the coefficients $\psi_l(u_j)\psi_m(u_j)\psi_m(u_j)\psi_o(u_j)$ are replaced by $(\hat \psi_{l,p}-\psi_l)\hat \psi_{m,p}\hat \psi_{n,p}\hat \psi_{o,p}$ [note that $A_{N,T}^*=B_{N,T}^*=0$ since the coefficients of the bootstrap process do not possess any time dependence]. If we employ \eqref{unglboots} and combine it with the fact that $|\hat {\underline d}- \underline d|<\alpha/4$ on the set $A_T(\alpha)$, we get
\begin{align}
\label{2.1abootstrap}
|\hat \psi_{l,p}-\psi_{l} | \leq C {p^4 \log(T)^{3/2}}T^{-1/2}|l|^{\alpha/4+ \underline d - 1} \quad \forall l \in \N ,
\end{align}

which together with \eqref{apprpsi} and the assumptions of the theorem implies
\begin{align}
\label{2.1abootstrapb}
|\hat \psi_{l,p} | \leq C |l|^{\alpha/4+ \underline d - 1} \quad \forall l,p \in \N .
\end{align}

Note that the coefficients in the MA($\infty$) representation of the bootstrap processes do not depend on time and that for such processes we only required 
\be \label{wesentlungl}
\sup_u|\psi_l(u)|\leq C|l|^{\underline d-1}
\ee
in the proof of Theorem \ref{thmapp1} a) to obtain appropriate bounds for the error. By using \eqref{2.1abootstrap} and \eqref{2.1abootstrapb} instead of \eqref{wesentlungl} and similar arguments as given in the proof of Theorem \ref{thmapp1} a)  it follows that
\bea
\E((\hat F_{1,T}^{*}-\hat F_{1,T,a}^{*})1_{A_T(\alpha)} | X_{1,T},...,X_{T,T})=F_{1,T}^{*,-}+O\Big(N^{4 \underline d-1+\alpha}p^4\log(T)^{3/2}T^{-1/2}\Big),
\eea

where $F_{1,T}^{*,-}$ is defined  as $F_{1,T}$ but with $f(u,\lambda)$ replaced by
\bea
\frac{\sigma_p^2}{2\pi}\sum_{l,m,n,o=-\infty}^\infty \hat \psi_{l,m,n,o,p} \exp(-i \lambda (l-m+n-o)) \times 1_{A_T(\alpha)}.
\eea

Since $\hat F_{2,T}^{*}$ and $\hat F_{2,T,2}^{*}$ are treated analogously, the claim follows [note that $F_{1,T}^{*,-}$ cancels out since the coefficients do not possess any time dependence]. $\hfill \Box$ \\

{\bf Proof of part b):} The assertion follows by similar arguments as given in the proof of Theorem \ref{thmapp1} b) employing \eqref{2.1abootstrap} and \eqref{2.1abootstrapb} instead of \eqref{wesentlungl} as above. The details are omitted for the sake of brevity.$\hfill \Box$ \\

\subsection{Proofs of the results in Section 3 and 4:}
{\bf Proof of Theorem \ref{thm1}:} The claim follows by employing the Cram\'{e}r-Wold device in combination with Theorem \ref{thmapp1}. $\hfill \Box$ \\

{\bf Proof of Theorem \ref{thm3}:}
Similarly to the proof of Theorem \ref{thmapp1}, the  two equations
\begin{eqnarray*}
\E (\hat{\tau}^2_1)&=&\frac{1}{\pi M}\sum_{j=1}^M\int_{-\pi}^{\pi}  f^4(u_j,\lambda)+ O \Big(\frac{1}{N^{1-8d_{\infty}}} \Big) ,\\
 \var(\hat{\tau}^2_1) &=& O\Big( \frac{1}{MN^{1-8d_{\infty}}} \Big)
\end{eqnarray*}
can be established. By Markov's inequality the assertion of the theorem follows. $\hfill \Box$ \\

{\bf Proof of Theorem \ref{thmstat2}:}  We define $\hat D_{T,a}^{2,*}$ as $\hat D_{T,2}^{2,*}$ and $\hat D_{T,a}^2$ as $\hat D_T^2$ but with $X_{t,T}$ replaced by $X_t(t/T)$ from \eqref{localappr}. Then part a) is obvious, because we have $\psi_l=\psi_l(u)$ for all $u \in [0,1]$ under the null hypothesis and $Z_t$ and $Z_t^*$ are both independent and standard normal distributed. Part b) follows from the proof of Theorem \ref{thmapp1}, so we focus on part c) and d). Note that \eqref{apprpsifolgerung} and Theorem \ref{thmapp1} a), b) imply
\bea
C_1N^{\max(8 \underline d -1,0)}/T \leq \Var(\hat D_{T,a}^*)/T \leq C_2(N^{\max(8 \underline d -1,0)}+\log(N)1_{\{\underline d = 1/8\}})/T
\eea

which directly yields part d). If we have
\be \label{patasymptotic}
P(A_T(\alpha))\rightarrow 1 \quad \text{ as } T \rightarrow \infty
\ee
for every $\alpha>0$, Part c) follows from Theorem \ref{thmapp2}, \eqref{assNM2}, the conditions on the rate of $p(T)$ if $\alpha$ is chosen sufficiently small. Finally, \eqref{patasymptotic} is a consequence of Lemma 4.3 in \cite{prevet2012}.  $\hfill \Box$ \\

{\bf Proof of Theorem \ref{mallowmetric}:} By employing the triangle inequality we can bound the Mallow metric between $\hat D_T^2/\sqrt{\Var(\hat D_T^2)}$ and $\hat D_T^{2,*}/\sqrt{\Var(\hat D_T^{2,*})}$ by
\bea
&&d_2 \Big(\hat D_T^2/\sqrt{\Var(\hat D_T^2)},\hat D_{T,a}^{2}/\sqrt{\Var(\hat D_{T,a}^{2})}\Bigr)+d_2 \Big(\hat D_{T,a}^{2}/\sqrt{\Var(\hat D_{T,a}^{2})},\hat D_{T,a}^{2,*}/\sqrt{\Var(\hat D_{T,a}^{2,*})}\Bigr)\\
&&+d_2 \Big(\hat D_{T,a}^{2,*}/\sqrt{\Var(\hat D_{T,a}^{2,*})},\hat D_T^{2,*}/\sqrt{\Var(\hat D_T^{2,*})}\Bigr)
\eea

[where $\hat D_{T,a}^2$ and $\hat D_{T,a}^{2,*}$ are the random variables from Theorem \ref{thmstat2} specified in the proof of which]. It follows from the proof of Theorem \ref{thmapp1} that the first summand converges to zero and the second summand equals zero because of  Theorem \ref{thmstat2} a). So it suffices to treat the third summand which is bounded by 
\bea
2\E \Big(\hat D_{T,a}^{2,*}/\sqrt{\Var(\hat D_{T,a}^{2,*})}- \hat D_T^{2,*}/\sqrt{\Var(\hat D_{T,a}^{2,*})}\Big)^2+2\E \Big(\hat D_{T}^{2,*}/\sqrt{\Var(\hat D_{T,a}^{2,*})}- \hat D_T^{2,*}/\sqrt{\Var(\hat D_T^{2,*})}\Big)^2.
\eea

We obtain from Theorem \ref{thmapp1} a) and b) that   a constant $L>0$ exists such that 
\be \label{varianzlowerbound}
\Var(\hat D_{T,a}^{2,*}) \geq LN^{8d_\infty -1}T^{-1}
\ee
[note the we are under the null hypothesis and that we therefore have  $d_\infty = \underline d$]. This combined with Theorem  \ref{thmapp2}, \eqref{patasymptotic} and the conditions on the growth rate on $p=p(T)$ yields that we can restrict ourselves to the second term, which is [up to the constant $2$] bounded by
\bea
&&\frac{\E(( \hat D_T^{2,*})^2)}{\Var(\hat D_{T,a}^{2,*})\Var(\hat D_{T}^{2,*})}\Big(\sqrt{\Var(\hat D_{T,a}^{2,*})}- \sqrt{\Var(\hat D_T^{2,*})}\Big)^2 \\
&\leq & \frac{\E(( \hat D_T^{2,*})^2)}{\Var(\hat D_{T,a}^{2,*})\Var(\hat D_{T}^{2,*})} \Big| \Var(\hat D_{T,a}^{2,*})- \Var(\hat D_T^{2,*}) \Big|.
\eea

If we follow the proof of Theorem \ref{thmapp1} a), b) and employ \eqref{2.1abootstrapb} and \eqref{patasymptotic}, we obtain that  
\be \label{abschzweitemoment}
\E(( \hat D_T^{2,*})^2)=O(\log(N)N^{\max(8 \underline d -1,0)+2\alpha}T^{-1}+N^{8 \underline d +2\alpha -2})
\ee
holds for every fixed $\alpha>0$. By employing \eqref{2.1abootstrap} and similar arguments as in the proof of Theorem \ref{thmapp2} we obtain thereafter
\bea
\Big| \Var(\hat D_{T,a}^{2,*}1_{A_T(\alpha)})- \Var(\hat D_T^{2,*}1_{A_T(\alpha)}) \Big|=O(p^8\log(T)^3\log(N)^2(N^{\max(8 \underline d-1,0)+2\alpha}T^{-2}+N^{8\underline d -2+2\alpha}T^{-1})).
\eea

The assertion then follows with \eqref{patasymptotic}--\eqref{abschzweitemoment} and the assumptions on the growth rate of $p=p(T)$. $\hfill \Box$

\section{Appendix: Auxiliary Lemmas}
\def\theequation{7.\arabic{equation}}
\setcounter{equation}{0}

Finally we show some lemmas which were employed  in the above proofs.

\begin{lem}
\label{Lemma 7.1}
Suppose $\mu, \nu, a, b \in \R$. Then there exists a constant $C \in \R$ such that the following holds:
\begin{itemize}
\item[a)] If $\mu , \nu >0$ and $b > a$, then
\begin{eqnarray}
\label{Gleichung 3.196(3)}
\sum_{\substack{k=0 \\k-a \geq 1 \\ -k + b \geq 1 }}^{N-1} \frac{1}{(k-a)^{1-\mu}} \frac{1}{(b-k)^{1-\nu}}
\leq \sum_{k=1+a}^{b-1} \frac{1}{(k-a)^{1-\mu}} \frac{1}{(b-k)^{1-\nu}}\leq  \frac{C}{(b-a)^{1-\mu-\nu}}.
\end{eqnarray}
\item[b)] If $0 <\mu,\nu$ and $0 < 1-\mu-\nu $, then it follows for $|a+b| > 0$
\begin{eqnarray}
\sum_{\substack{k=1 \\ k+b \geq 1 \\ k-a \geq 1}}^{N-1} \frac{1}{(k+b)^{1-\mu}} \frac{1}{(k-a)^{1-\nu}} \leq 
\sum_{\substack{k=1 \\ k+b \geq 1 \\ k-a \geq 1}}^{\infty} \frac{1}{(k+b)^{1-\mu}} \frac{1}{(k-a)^{1-\nu}} 
\leq \frac{C}{|a+b|^{1-\mu-\nu}}.
\label{eq1}
\end{eqnarray}
\item[c)] If $0 < \nu < 1-\mu$ and $y,z \geq 1$, then
\begin{eqnarray}
\label{logarithmus und 1+y}
&&\sum_{k=1+y}^{\infty} \frac{\log (k)}{k^{1-\mu}} \frac{1}{(k-y)^{1-\nu}}\leq  C\frac{\log(y)}{y^{1-\mu-\nu}}, \\
\label{logarithmus und 1}
&&\sum_{k=1}^{\infty} \frac{\log(k+z)}{(k+z)^{1-\mu}} \frac{1}{k^{1-\nu}} \leq C \frac{\log(z)}{z^{1-\mu-\nu}}.
\end{eqnarray}
\end{itemize}
\end{lem}
{\bf Proof:} a) Using equation 3.196(3) in \cite{grad1980}, it follows that
\begin{eqnarray*}
\sum_{k=1+a}^{b-1} \frac{1}{(k-a)^{1-\mu}} \frac{1}{(b-k)^{1-\nu}} \leq  \int_{a}^{b} \frac{1}{(x-a)^{1-\mu}}\frac{1}{(b-x)^{1-\nu}} \,dx \lesssim \frac{1}{(b-a)^{1-\mu-\nu}}.
\end{eqnarray*}

b) If $a+ b >0$ we can bound the middle sum in \eqref{eq1} by
\begin{eqnarray*}
&&
\sum_{\substack{k= \max\{ 1,1-b,1+a \}}}^{\infty} \frac{1}{(k+b)^{1-\mu}} \frac{1}{(k-a)^{1-\nu}}\leq \sum_{\substack{k=1+a }}^{\infty} \frac{1}{(k+b)^{1-\mu}} \frac{1}{(k-a)^{1-\nu}}\\
& \leq & \int_{a}^{\infty} \frac{1}{(x+b)^{1-\mu}} \frac{1}{(x-a)^{1-\nu}}\,dx \lesssim \frac{1}{(a+b)^{1-\mu-\nu}}.
\end{eqnarray*}
The last inequality follows from the equations 3.196(2) and 3.191(2) [for choosing $b=0$] in \cite{grad1980}.
On the other hand, if $a+b < 0$ we have
\begin{eqnarray*}
&& \sum_{\substack{k= \max\{ 1,1-b,1+a \}}}^{\infty} \frac{1}{(k+b)^{1-\mu}} \frac{1}{(k-a)^{1-\nu}} \leq  \sum_{\substack{k=1-b}}^{\infty} \frac{1}{(k+b)^{1-\mu}} \frac{1}{(k-a)^{1-\nu}}\\
&=&   \sum_{\substack{k=1+(-b)}}^{\infty} \frac{1}{(k-(-b))^{1-\mu}} \frac{1}{(k+(-a))^{1-\nu}} \lesssim  \frac{1}{(-a-b)^{1-\mu-\nu}}.
\end{eqnarray*}
The last inequality follows with \cite{grad1980} as above.\\

c) We start with (\ref{logarithmus und 1+y}).
Using equation 13.2(18) in \cite{erdelyi2} yields
\begin{eqnarray*}
\sum_{k=1+y}^{\infty} \frac{\log (k)}{k^{1-\mu}} \frac{1}{(k-y)^{1-\nu}} \leq \int_{y}^{\infty} \frac{\log (x)}{x^{1-\mu}} \frac{1}{(x-y)^{1-\nu}} \,dx \lesssim \frac{\log(y)}{y^{1-\mu-\nu}}.
\end{eqnarray*}
Concerning (\ref{logarithmus und 1}) we use equation 6.4(23) in \cite{erdelyi1} which implies 
\begin{eqnarray*}
\sum_{k=1}^{\infty} \frac{\log(k+z)}{(k+z)^{1-\mu}} \frac{1}{k^{1-\nu}} \leq \int_{0}^{\infty} \frac{\log(x+z)}{(x+z)^{1-\mu}} \frac{1}{x^{1-\nu}}\,dx \lesssim  \frac{\log(z)}{z^{1-\mu-\nu}}.
\end{eqnarray*}
$\hfill \Box$

%
%

\begin{lem}
\label{j_1 ungleich j_2}
If $0 < d_{\infty} < \frac{1}{4}$, then
\begin{eqnarray*}
&&\frac{1}{M^2N^4} \sum_{\substack{j_1,j_2=1 \\ j_1 \not = j_2}}^M \sum_{\substack{r_1,s_1=0  }}^{N-1} \sum_{\substack{v_1,w_1,x_1,y_1=0  \\ v_1-w_1+r_1-s_1=0 \\ 0 \leq |v_1-w_1| \leq N-1 \\ }}^{\infty} \sum_{\substack{v_2,w_2,x_2,y_2 =0\\ r_1+ v_2-x_1-s_1-w_2+y_1 -y_2+x_2 =0\\ 0 \leq r_1+ v_2-x_1+({j_1}-{j_2})N \leq N-1\\ 0 \leq  s_1+w_2-y_1+(j_1-j_2)N \leq N-1 \\ 0 \leq |y_2-x_2| \leq N-1}}^{\infty} \\ 
&& \psi_{v_1}(u_{j_1})\psi_{w_1}(u_{j_1})\psi_{x_1}(u_{j_1})\psi_{y_1}(u_{j_1})\psi_{v_2}(u_{j_2})\psi_{w_2}(u_{j_2}) \psi_{x_2}(u_{j_2})\psi_{y_2}(u_{j_2}) \\
&&(N^2-N|y_2-x_2|-N|v_1-w_1|+|v_1-w_1|y_2-x_2| )  = O \Big(\frac{\log(N)}{N^{1-8d_{\infty}}T}  \Big).
\end{eqnarray*}
\end{lem}
{\bf Proof:}
Firstly, we set $0 \leq w_1= r_1-s_1+v_1$ and $0 \leq x_2= s_1-r_1-v_2+x_1+w_2-y_1+y_2$. Then we define $p:= r_1+ v_2-x_1+({j_1}-{j_2})N$ and rearrange to $ 0 \leq x_1= r_1-p+v_2+({j_1}-{j_2})N$.
Since $p \in \{ 0,1,2 \ldots, N-1\}$, it follows that if $p,r_1,v_2,x_1,j_1$ are fixed, there are at most two possible values for $j_2$. Hence it is enough to consider the following expression with $ 1 \leq |C_1| \leq M-1$
\begin{eqnarray*}
&&\frac{1}{MN^4}  \sum_{\substack{r_1,s_1=0  }}^{N-1} \sum_{\substack{v_1,w_1,x_1,y_1=0  \\ v_1-w_1+r_1-s_1=0 \\ 0 \leq |v_1-w_1| \leq N-1 \\ }}^{\infty} \sum_{\substack{v_2,w_2,x_2,y_2 =0\\ r_1+ v_2-x_1-s_1-w_2+y_1 -y_2+x_2 =0\\ 0 \leq r_1+ v_2-x_1+C_1N \leq N-1\\ 0 \leq  s_1+w_2-y_1+C_1N \leq N-1 \\ 0 \leq |y_2-x_2| \leq N-1}}^{\infty} \\ 
&& \psi_{v_1}(u_{j_1})\psi_{w_1}(u_{j_1})\psi_{x_1}(u_{j_1})\psi_{y_1}(u_{j_1})\psi_{v_2}(u_{j_2})\psi_{w_2}(u_{j_2}) \psi_{x_2}(u_{j_2})\psi_{y_2}(u_{j_2}) \\
&&(N^2-N|y_2-x_2|-N|v_1-w_1|+|v_1-w_1|y_2-x_2| ) \\
&\stackrel{\eqref{apprpsi}}{\lesssim}&\frac{1}{MN^2} \sum_{\substack{p,r_1,s_1=0  }}^{N-1} 
\sum_{\substack{v_1,x_1,y_1=1  \\ 1 \leq r_1-s_1+v_1 }}^{\infty} \sum_{\substack{v_2,w_2,x_2,y_2 =1 \\ 1 \leq r_1-p+ v_2+C_1N\\  1 \leq s_1-p+w_2+y_2-y_1+C_1N}}^{\infty} \frac{1}{v_1^{1-d_{\infty}}}\frac{1}{y_1^{1-d_{\infty}}}
\frac{1}{v_2^{1-d_{\infty}}}\frac{1}{w_2^{1-d_{\infty}}}\frac{1}{y_2^{1-d_{\infty}}} \\ 
&&\frac{1}{(r_1-s_1+v_1)^{1-d_{\infty}}}\frac{1}{(r_1-p+ v_2+C_1N)^{1-d_{\infty}}} \frac{1}{(s_1-p+w_2+y_2-y_1+C_1N)^{1-d_{\infty}}}\\
&\stackrel{(\ref{Gleichung 3.196(3)})}{\lesssim}&\frac{1}{MN^2} \sum_{\substack{p,r_1=0  }}^{N-1} 
\sum_{\substack{v_1,x_1,y_1=1   }}^{\infty} \sum_{\substack{v_2,w_2,x_2,y_2 =1\\  1 \leq r_1-p+ v_2+C_1N \\ 1 \leq r_1-p+w_2-y_1+y_2+v_1+C_1N}}^{\infty} \frac{1}{v_1^{1-d_{\infty}}}\frac{1}{y_1^{1-d_{\infty}}}
\frac{1}{v_2^{1-d_{\infty}}}\frac{1}{w_2^{1-d_{\infty}}} \frac{1}{y_2^{1-d_{\infty}}} \\ 
&&\frac{1}{(r_1-p+ v_2+C_1N)^{1-d_{\infty}}}\frac{1}{(r_1-p+w_2-y_1+y_2+v_1+C_1N)^{1-2d_{\infty}}}\\
&\stackrel{(\ref{Gleichung 3.196(3)}), (\ref{eq1})}{\lesssim}&\frac{1}{MN^2}  \sum_{\substack{p,r_1=0  }}^{N-1} \frac{1}{|r_1-p+C_1N|^{2-8d_{\infty}}} \lesssim \frac{1}{MN^2}\sum_{\substack{p,r_1=0  }}^{N-1} 
\frac{1}{(r_1-p+N)^{2-8d_{\infty}}} \\
&\lesssim &\frac{\log(N)}{MN^{2-8d_{\infty}}} =\frac{\log(N)}{N^{1-8d_{\infty}}T}.
\end{eqnarray*}
$\hfill \Box$

\end{document}